\documentclass[12pt]{amsart}
\usepackage{amscd,amsthm,amsfonts,amssymb,amsmath,euscript}
\usepackage{marvosym}
\usepackage[dvips]{graphicx}
\usepackage[dvips]{graphics}
\usepackage[matrix,arrow]{xy}
\usepackage{longtable}


\newtheorem{lemma}{Lemma}[section]
\newtheorem{theorem}[lemma]{Theorem}

\newtheorem{proposition}[lemma]{Proposition}
\newtheorem{conjecture}[lemma]{Conjecture}

\theoremstyle{definition}

\newtheorem{definition}[lemma]{Definition}

\theoremstyle{remark}
\newtheorem{remark}[lemma]{Remark}

\newcommand{\QQ}{{\mathbb Q}}
\newcommand{\ZZ}{{\mathbb Z}}

\newcommand{\PP}{{\mathbb P}}
\newcommand{\CC}{{\mathbb C}}
\newcommand{\RR}{{\mathbb R}}

\newcommand{\Aff}{{\mathbb A}}

\newcommand{\Fuk}{{\operatorname{{\sf Fuk}}}}

\newfont{\smallskob}{cmbx7 scaled\magstep4}
\newfont{\bigskob}{cmbx12 scaled\magstep4}

\newcommand{\Spec}{\mathrm{Spec}\,}
\newcommand{\Gap}{\mathrm{Gap}\,}
\newcommand{\NLGap}{\mathrm{NLGap}\,}

\newcommand{\cO}{\mathcal{O}}

\def\PP{{\mathbb P}}
\def\CC{{\mathbb C}}
\def\PP{{\mathbb P}}
\def\QQ{{\mathbb Q}}
\def\ZZ{{\mathbb Z}}

\def\cO{\mathcal{O}}

\def\cO{\mathcal{O}}
\def\cP{\mathcal{P}}

\newlength{\nanowidth} 

\title{Double solids, categories and non-rationality}

\author{Atanas Iliev,  Ludmil Katzarkov, Victor Przyjalkowski}

\begin{document}

\maketitle

\begin{flushright}
To Slava  --- teacher and friend --- with admiration.
\end{flushright}

\begin{abstract}

This paper suggests a new approach to questions of rationality
of threefolds  based on  category theory.
Following \cite{BFK} and \cite{FK} we enhance constructions from \cite{KUZ} by introducing Noether--Lefschetz spectra --- an interplay between Orlov spectra \cite{O} and Hochschild homology.
The main goal  of this  paper is to suggest a series of interesting examples where  above techniques might apply.
 We start by constructing  a sextic double solid $X$ with
$35$ nodes and torsion in $H^3(X,\ZZ)$. This is a novelty --- after the classical example of Artin
and Mumford (1972), this is the second example of a Fano threefold with a torsion in  the 3-rd integer homology group. In particular $X$ is non-rational.
We consider other examples as well ---
 $V_{10}$ with  10 singular points and
 double covering of quadric ramified in octic with 20 nodal singular points.

After analyzing the geometry of their Landau Ginzburg models  we suggest a general non-rationality picture based on Homological Mirror Symmetry and category theory.

\end{abstract}

\section{ Introduction}

This paper suggests a new approach to questions of rationality
of threefolds  based on  category theory.
It was inspired by recent work of V.\,Shokurov and by A.\,Kuznetsov's idea about the Griffiths component (see~\cite{Kuz08}). This work is a natural continuation of ideas developed in \cite{KRIS1},
\cite{GAR} and of ideas of Kawamata and his school.

We first  extend  classical example of Artin and Mumford
to construct a sextic double solid $X$ with
$35$ nodes and torsion in $H^3(X,\ZZ)$. The construction is based on
an approach by M.\,Gross and suggests close relation between   Artin and Mumford  example and the sextic double solid $X$ with
$35$ nodes. This example, a novelty on its own,  opens a possibility of series of interesting examples ---
$V_{10}$ with  10 singular points and
double covering of quadric ramified in octic with 20 nodal singular points.

 In this paper  we start investigating these examples  from the point of view of Homological Mirror Symmetry (HMS). We consider the mirrors  of the
 sextic double solid $X$ with
$35$ nodes, of the  Fano variety $V_{10}$ with 10 singular points in general position and of the
 double covering of quadric ramified in octic with 20 nodal singular points.
  We note that the monodromy around the singular fiber over zero of the  Landau--Ginzburg models is strictly unipotent in all these examples, which
suggests that the
categorical behavior should be very similar to the one of the  Artin--Mumford example.
We conjecture that the reason for   categorical similarity in all these examples is that they contained the category of an   Enriques surface as a semiorthogonal summand in their derived categories.  This is done  in Section 5,  where we introduce Landau--Ginzburg models and compare their singularities.

In Section 6 we introduce several new rationality invariants coming out of the
notions of spectra and enhanced Noether--Lefschetz  spectra of categories.  We give a conjectural categorical explanation of the examples from Sections 2, 3, 4, 5. The novelty (conjecturally)  is that non-rationality of these examples cannot be picked by Orlov spectra but it is detected by the Noether--Lefschetz spectrum.

The paper is organized as follows:  in Sections 2, 3, 4  we describe  classical calculations of  a sextic double solid $X$ with
$35$ nodes. Section 5 contains some  mirror considerations studying some Landau--Ginzburg models.
Section 6 suggests a  general categorical framework for studying the phenomena in Sections 2--5.

The paper is based on examples we have analyzed in
\cite{MARNON}, \cite{KNS}, \cite{FK1}, \cite{FK}, \cite{KK}. All these  suggest a direct connection between monodromy of Landau--Ginzburg models, spectra and wall crossings  in the moduli space of stability conditions, which was partially explored in \cite{IKS}.  This paper is a humble attempt to shed some light on this connection. We expect that further application of this method will be the theory of three dimensional conic bundles --- a small part of  huge algebro--geometric heritage  of V.\,Shokurov (see Remark~\ref{remark:quadric solid}). In particular we expect that
Noether--Lefschetz  spectra of categories would allow us to prove nonrationality of new classes of  conic bundles --- classes where the method of Intermediate Jacobian does not work.

\medskip

{\bf Acknowledgements.}

The authors thank  D.\,Favero, G.\,Kerr, M.\,Kontsevich, A.\,Kuznetsov, D.\,Orlov, T.\,Pantev for useful discussions,
the referee for proofreading the paper,
 and C.\,Shramov for pointing out the example of double covering of quartic ramified in octic with 20 nodal points.
A.\,I.  was funded by  FWF  Grant P20778.
L.\,K. was funded by NSF Grant DMS0600800, NSF FRG Grant DMS-0652633, FWF
Grant P20778 and an ERC Grant --- GEMIS.
V.\,P. was funded by grants NSF FRG DMS-0854977, NSF DMS-0854977, and NSF DMS-0901330; grants FWF P 24572-N25 and FWF P20778; RFFI grants 11-01-00336-a, 11-01-00185-a, 12-01-33024, and 12-01-31012;
grants MK-1192.2012.1, NSh-5139.2012.1, and AG Laboratory GU-HSE RF government grant, ag. 11 11.G34.31.0023.

\bigskip

All varieties   considered in this paper  are defined   over the  field of complex numbers   $\CC$. The torsion subgroup of given  group $G$ is denoted by $Tors\,(G)$; the $n$-torsion subgroup is denoted by $Tors_n\,(G)$.
We denote du Val singularities of ADE type by $A_n$, $D_n$, and $E_n$. We denote a Landau--Ginzburg model of  a variety $X$ by $LG(X)$.

\section{Determinantal double solids and Brauer--Severi varieties}\label{ai}

\subsection{The classical Artin--Mumford example}\label{ai-1}
A double solid is an irreducible double covering $\pi: X \rightarrow \PP^3$.
The branch locus of such $\pi$ is a surface $S \subset \PP^3$ of even degree.
In 1972 Artin and Mumford gave an example of a special
singular quartic double solid $X$ (i.e. $deg \ S = 4$) which is non-rational
because of the existence of a non-zero
2-torsion in its integer cohomology group $H^3(X,\ZZ)$, see \cite{AM}.
Since quartic double solids are unirational (see, for instance,~\cite{IP99}, Example~10.1.3(iii)), this gives
(together with the examples presented at the same time by Iskovskikh--Manin
and Clemens--Griffiths) an example of a non-rational unirational threefold.

In \cite{AMG} Aspinwall, Morrison and Gross present a special case of a
singular Calabi--Yau threefold --- an octic double solid $X$
(i.e. $deg\ S = 8$) with $80$ nodes on $S$ and a non-zero 2-torsion
in $H^3(X,\ZZ)$.

In this section  we adapt an  approach used in \cite{AMG} to check again
the existence of the 2-torsion in $H^3(X,\ZZ)$
for the Artin--Mumford quartic double solid $X$, and present an example of a sextic double solid
$X$ with 35 nodes and a non-zero torsion in $H^3(X,\ZZ)$.
In particular this special nodal sextic double solid
is not rational. Other examples are presented in sections to follow.

\subsection{Quadric bundles and determinantal double solids}\label{ai-2}
Let $X_0$ be a smooth complex projective variety,
let $L$ be an invertible sheaf on $X_0$,
and let $E \rightarrow X_0$ be a vector bundle of rank $r \ge 2$
over $X_0$.

A {\it quadric bundle} in $E$ parameterized by $L$ is an $\cO_{X_0}$-map
$$
\varphi\colon L^{-1} \rightarrow Sym^2 E^*.
$$
The determinantal loci of $\varphi$ are  subvarieties
$$
D_{r-k} = D_{r-k}(\varphi) = \{ x \in X_0 : rank\ \varphi_x \le r-k \},\  k = 0,1,2,\ldots
$$
Geometrically a quadric bundle $\varphi$ represents a bundle of quadrics
$${\mathcal Q} = \{ Q_x \subset \PP(E_x) : x \in X_0 \},$$
and
$$
D_{r-k} = \{ x \in X_0 : rank\ Q_x \le r-k \}.
$$
If $D_{r-k} \subset X_0$ are nonempty and have the expected codimensions $k(k+1)/2$ then
their classes in $A_*(X_0)$ can be computed by the formulas in
\cite{HT} or \cite{JLP}.
For our purposes we need only to know  explicit formulas for
first two determinantals $D_{r-1}$ and $D_{r-2}$,
which can be computed formally as follows.
Rewrite $\varphi$ in the form
$$
\varphi\colon \cO_{X_0} \rightarrow Sym^2 (E^* \otimes L^{1/2}),
$$
and compute $c(E^* \otimes L^{1/2}) = 1 + c_1 + c_2 + \ldots + c_r$.
Then
\begin{equation*}\label{determinantals}
D_{r-1} = 2c_1 \ \mbox{ and } \ D_{r-2} = 4(c_1c_2-c_3).
\end{equation*}

In  particular case when the base $X_0 = \PP^n$ is a projective space,
then the determinantal locus $D_{r-1}$ is a hypersurface in
$\PP^n$ of even degree; therefore $D_{r-1}$ defines a double
covering
$$
\pi\colon X \rightarrow \PP^n
$$
branched along $D_{r-1}$. We call such $X$ a {\it determinantal double solid}.

\subsection{Cohomological Brauer groups and Brauer--Severi varieties}\label{ai-3}
Let $X$ be a complex algebraic variety, let $\cO_X$ be the structure
sheaf of $X$, and let $\cO_X^{*}$ be the sheaf of units in $\cO_X$.
The Picard group and the (cohomological) Brauer group of $X$
are correspondingly the 1-st and the 2-nd cohomology groups
$$
Pic(X) = H^1(X,\cO^*_X) \ \mbox{ and } \ Br(X)=H^2(X,\cO^*_X).
$$
There is an exact sequence
$$
Pic(X) \otimes \QQ/\ZZ \rightarrow H^2(X,\QQ/\ZZ) \rightarrow Br(X) \rightarrow 0,
$$
see 3.1 in part II of \cite{Gr}.
If in addition $X$ is non-singular and it fulfills  conditions
\begin{equation}\label{h3-conditions}
Pic(X) = H^2(X,\ZZ) \ \mbox{ and } \ H^1(X, \cO_X) = H^2(X,\cO_X) = 0,\
\end{equation}
then by the universal coefficient theorem
$Br(X) \cong Tors\,(H^3(X,\ZZ))$, see e.g. \cite{AM}.
For any $X$ as above,
{\it a Brauer--Severi variety} over $X$ is
a variety ${\cP}$ with a structure of a $\PP^n$-bundle
$f\colon{\cP} \rightarrow X$ over $X$.

Not any Brauer--Severi variety
is a projectivisation of a vector bundle over $X$, and the
Brauer group gives obstructions for a Brauer--Severi variety
to be a presented as a projectivisation of such.
On $X$, we consider   exact sequence
$$
0 \rightarrow \cO^*_X \rightarrow GL_{n+1} \rightarrow PGL_{n+1} \rightarrow 0,
$$
where $\cO^*_X$ is the multiplicative group of $X$.

The corresponding  long exact sequence is
$$
0 \longrightarrow Pic(X) \longrightarrow
 H^1(X,GL_{n+1})
\stackrel{j}\longrightarrow H^1(X,PGL_{n+1})
\stackrel{\delta}\longrightarrow  Br(X)
\longrightarrow \ldots
$$
The vector bundles $E \rightarrow X$ of rank $(n+1)$ are  elements of
the cohomology group $H^1(X,GL_{n+1})$, while the $\PP^n$-bundles
${\cP} \rightarrow X$ are  elements  of $H^1(X,PGL_{n+1})$.

Therefore by  above sequence
%
the $\PP^n$-bundle ${\cP}$ is not a projectivisation of a vector bundle
on $X$ iff $\delta({\cP}) \not= 0$.
 Since $(n+1)\delta = 0$ then
any ${\cP}$ with $\delta(\cP) \not= 0$ gives rise to a non-zero
(n+1)-torsion element $\delta({\cP}) \in Br(X)$.
If moreover $X$ fulfills  conditions (\ref{h3-conditions})
then ${\cP}$ represents a non-zero (n+1)-torsion element of $H^3(X,\ZZ) \cong Br(X)$.
In  particular case we consider  below ${\cP}$ is a $\PP^1$-bundle
which is not a projectivisation of a vector bundle, thus representing
a non-zero 2-torsion element of $H^3(X,\ZZ)$.

In the next sections we will use the following:

\begin{lemma}\label{torsion-criterion} {\bf Torsion criterion for non-rationality}.
For the smooth complex variety $Y$ the torsion subgroup $Tors\,(H^3(Y,\ZZ))$
is a birational invariant of $X$.
In particular if $Y$ is rational then $Tors\,(H^3(Y,\ZZ)) = 0$.
\end{lemma}

\begin{proof}
See Proposition~1 in \cite{AM} or \S 9 in \cite{Be}.
\end{proof}

\section{Determinantal sextic double solid $X$ with a non-zero 2-torsion in
         $H^3(X,\ZZ)$}\label{ai-4}

\subsection{The double solids of Artin--Mumford, Aspinwall--Morrison--Gross,
and  determinantal sextic double solid}
The Artin--Mumford threefold from \cite{AM}
is a special double solid with a branch locus
--- a quartic surface $S$ with 10 nodes and with a torsion in the
3-rd integer cohomology group $H^3 = H^3(\widetilde{X},\ZZ)$,
where $\widetilde{X} \rightarrow X$ is the blowup of $X$ at its nodes.
As it was shown later by Endrass, the group $H^3$
of a double solid $X$ branched over a nodal quartic surface $S$
can have a non-zero torsion only in case when $S$ has $10$ nodes,
see~\cite{En}. Therefore the branch loci of eventual further examples of nodal 3-fold
double solids with a non-zero torsion in the 3-rd integer cohomology group $H^3$ should be of degree $d$ either equals $2$ or $\ge 6$. If in addition we require such $X$
to be a Fano threefold then $d$ must be $\le 6$, i.e. if exists such
$X$ must be a sextic double solid or a double quadric.
Notice that non-singular Fano threefolds $X$ have a zero torsion
in $H^3 = H^3(X,\ZZ)$, so the requirement $X$ to be singular
(and nodal --- for simplicity) is substantial.

In \cite{AMG} Aspinwall, Morrison, and Gross study a special
case of a Calabi--Yau threefold which is a double solid $X$
with a torsion in $H^3$ and with a branch locus $S$ of degree $8$
(an octic double solid).
The similarity  between the Artin--Mumford quartic double solid and the
octic double solid from \cite{AMG} is that they both are
{\it determinantal} double solids. Both these varieties $X$ are
singular --- in the Artin--Mumford case $X$ has 10 ordinary double points
(nodes) while the octic double solid from \cite{AMG} has 80 nodes.

\medskip

Below we describe an example of a determinantal nodal sextic double
solid $X$ with a torsion in $H^3$.
After the example of Artin and Mumford,
this is the 2-nd example of a (necessary) singular nodal Fano threefold
(see above)
with a torsion in the 3-rd integer cohomology group.
In particular our $X$ must be non-rational, see Lemma~\ref{torsion-criterion}.

It is shown by Iskovskikh (see \cite{Is}) that the general sextic double solid
is non-rational due to the small group $Bir(X)$ of birational automorphisms
of $X$. This argument  has been extended later by Cheltsov and Park proving the non-rationality of certain singular sextic double solids, see \cite{CP}.

From this point of view, the example studied below is a non-rational
sextic double solid $X$ with 35 ordinary double points.
According to  Cheltsov (private communication), the non-rationality of this $X$cannot be derived, at least for now, from the results of \cite{CP}.

The proof of the non-rationality of $X$ presented below follows
ideas from Appendix in \cite{AMG}.

\subsection{The determinantal sextic double solid}
Let $\PP^3 \times \PP^4 \subset \PP^{19}$
be  Segre variety of $\CC^*$-classes of non-zero $4 \times 5$ matrices,
and let
$$
W =  (\PP^3 \times \PP^4) \cap H \cap F
$$
be a general complete intersection of $\PP^3 \times \PP^4$
with a hyperplane $H = \PP^{18} \subset \PP^{19}$
and a divisor $F$ of bidegree (1,2).
Let $Z = (\PP^3 \times \PP^4) \cap H$,
and denote by $p_Z$ and $p_W$ the restrictions of the projection
$p\colon \PP^3 \times \PP^4 \rightarrow \PP^3$
to $Z$ and to $W$.
The projection $p_W$ defines a structure of a quadric bundle
$$
p_W\colon W \rightarrow \PP^3
$$
on $W$ with fibers --- quadrics
$Q_x = p_W^{-1}(x)$ in the 3-spaces
$$
\PP^3_x = p_Z^{-1}(x) = (x \times \PP^4) \cap H, \ x \in \PP^3.
$$


The $\PP^3$-bundle $p_Z\colon Z \rightarrow \PP^3$ is a projectivisation
of the rank 4 vector bundle $E$ on $\PP^3$ defined by vanishing
linear form $H$ defining a hyperplane section $h$ on
fibers of $p\colon \PP^3 \times \PP^4 \rightarrow \PP^3$:
$$
0 \longrightarrow E \longrightarrow \cO_{\PP^3}^{\oplus 5}
  \stackrel{h}\longrightarrow \cO_{\PP^3}(1) \longrightarrow 0;
$$
therefore $c(E^*) = 1 + h + h^2 + h^3$ in $A_*(\PP^3) = \CC[h]/(h^4)$.
Since $W$ is an intersection of $Z = \PP(E) \rightarrow \PP^3$
with a bidegree (1,2) divisor, then the bundle of quadrics
defining a  quadric bundle
$p_W\colon W \rightarrow \PP^3$ is given by the map
$$
\varphi\colon \cO_{\PP^3}(-1) \rightarrow S^2 E^*.
$$
So
$c(E^*(\frac{1}{2})) = 1 + c_1 + c_2 + c_3 = 1 +  3h + 4h^2 + \frac{13}{4}h^3$,
and hence
$$
[D_3(\varphi)]  = 2c_1 = 6h \ \mbox{ and } \ [D_2(\varphi)] = 4(c_1c_2-c_3) = 35.
$$
For a general choice of a bidegree $(1,2)$ divisor $F$,
the branch locus
$$
S = D_3(\varphi)
$$
is a sextic surface in $\PP^3$
with 35 nodes --- the 35 points of
$$
\delta = D_2(\varphi) =  \{ p_1 ,\ldots,p_{35} \}.
$$
Let
$$
\pi\colon X \rightarrow \PP^3
$$
be a double covering branched along  sextic surface $S = D_3$.
Since $Sing\,(S) = \delta$, and the points $p_i \in \delta$ are nodes of $S$,
then  sextic double solid $X$ has $35$ nodes --- the preimages of the
35 points $p_1,\ldots,p_{35}$ of $\delta$.

\begin{proposition}\label{sds}
Let
$W =  (\PP^3 \times \PP^4) \cap H \cap F$
be a general complete intersection of $\PP^3 \times \PP^4$
with a hyperplane and a divisor of bidegree (1,2),
Then:

(1) the degeneration locus $S = D_3$ of quadric fibration
$p_W\colon W \rightarrow \PP^3$, induced by the projection
$p\colon \PP^3 \times \PP^4 \rightarrow \PP^3$,
is a sextic surface with 35 nodes;

(2) Let $\pi\colon X \rightarrow \PP^3$  be a  double covering
branched along the sextic surface $S = D_3$.
Then the group $H^3(X,\ZZ)$ contains a non-zero 2-torsion element;
in particular $X$ is non-rational.
\end{proposition}



\subsection{Proof of Proposition \ref{sds}}
Part (1) follows from previous considerations. It remains to verify (2).
Following an approach from \cite{AMG},  we will  find bellow a non-zero
2-torsion element of $H^3(X,\ZZ)$, by representing it as a
Brauer--Severi variety over the smooth part of $X$.
Together with Lemma \ref{torsion-criterion}
this  completes the proof.

We consider   quadric bundle $p_W\colon W \rightarrow \PP^3 = \PP^3(x)$,
and restrict it over  open subset
$$
\PP^3_0 = \PP^3 - \delta.
$$

We define
$$
S_0 = S - \delta, \
X_0 = X - \delta_X
\ \mbox{ and } \ W_0 = W - \delta_W,
$$
where $\delta_X = \pi^{-1}(\delta)$ is  isomorphic preimage of
$\delta = \{ p_1,\ldots,p_{35} \}$ on $X$, and $\delta_W = p^{-1}(\delta)$ is the set of
35 rank 2 quadric surfaces $Q_i = p^{-1}(p_i)$, $i = 1,\ldots,35$.
Outside $\delta_W$, the projection $p$ restricts to a quadric bundle
$$
p_{W_0}\colon W_0 \rightarrow \PP^3_0
$$
with degeneration locus $S_0$.

Let
$\pi\colon X_0 \rightarrow \PP^3_0$
be  induced determinantal double covering branched along $S_0$.
As it follows from our  construction the fibers of the quadric bundle
$p_{W_0}\colon W_0 \rightarrow \PP^3_0$
are  quadrics $Q_x \subset \PP^3_x, x \in \PP^3 - \delta$.

Let ${\cP}$ be the family of lines $l \subset W_0$  in
 the quadrics $Q_x, x \in \PP^3 - \delta$, and let
$\cP_0 \subset {\cP}$ be the family of these lines
$l \in {\cP}$ which lie on quadrics $Q_x, x \in \PP^k_0 - \delta$.

Let us denote by
$$
f_P\colon {\cP} \rightarrow \PP^3
$$
 the map sending a line $l \subset Q_x$ to a  point $x \in \PP^3$,
and let us   denote  by
$f_P\colon \cP_0 \rightarrow \PP^3_0$
its restriction over $\PP^3_0$.

We also define
$$
\pi_0\colon X_0 \rightarrow \PP^3_0
$$
to be the restriction of the double covering $\pi\colon X \rightarrow \PP^3$
to $X_0 = X - \delta_X$.

For any point $x \in \PP^3_0 - S_0 = \PP^3 - S$
the quadric $Q_x \subset \PP^3_x$ is smooth, while for any
$x \in S_0 = S - \delta$ the quadric $Q_x$ is a quadratic
cone of rank 3 in $\PP^3_x$.

Then we have

$$
f_P^{-1}(x) \cong P^1 \vee P^1 \ \mbox{ for } x \in \PP^3_0 - S_0 = \PP^3 - S,
$$
$$
\mbox{ and } \ f_P^{-1}(x) \cong \PP^1 \ \mbox{ for } x \in S_0  = S - \delta.
$$

Since $S_0$ is also  branch locus
of the double covering
$\pi_0\colon X_0 \rightarrow \PP^3_0$,
we identify   points of $X_0$ with  generators of
quadrics $Q_x, x \in \PP^3_0$. Therefore the mapping
$\cP_0 \rightarrow \PP^3_0$ is represented as a composition
$$
\cP_0 \stackrel{f_0}\longrightarrow X_0 \stackrel{\pi_0}\longrightarrow \PP^3_0,
$$
where
$$
f_0\colon \cP_0 \rightarrow X_0
$$
is a $\PP^1$-fibration sending the sets of lines $l$ on the quadrics $Q_x$
to the generators of $Q_x$ containing $l$.
Let
$$
\widetilde{X} \rightarrow X
$$
be the  blowup of $X$ at
35 nodes of $X$ identified with  35 double points
$p_1,\ldots,p_{35}$ of the surface $S$.
Following \cite{AMG}
we  see that $\cP_0$ is not a projectivisation of
a vector bundle over $X_0$. This  yields that
the Brauer group $Br(\widetilde{X})$ has a
non-zero element of order two, representing a non-zero
2-torsion element in $H^3(\widetilde{X},\ZZ)$.

\medskip

Suppose that $f_0\colon \cP_0 \rightarrow X_0$  is a
projectivisation of a rank 2 vector bundle $E \rightarrow X_0$.
Up to a twist by a line bundle, we can always assume that $E$
has sections.
Next, any section of $E$  gives rise to a rational section
of $f_0\colon \cP_0 = \PP(E) \rightarrow X_0$. The following lemma
concludes the argument:

\begin{lemma}\label{not-a-projectivization}
The $\PP^1$-fibration $f_0\colon \cP_0 \rightarrow X_0$ has no rational sections.
In particular $\cP_0$ is not a projectivisation of a rank 2 vector bundle
on $X_0$.
\end{lemma}

\begin{proof}[Proof (see \cite{AMG} for more detail).]
Suppose that $f_0$ has a rational section, i.e. a rational map
$\sigma\colon X_0 \rightarrow \cP_0$ defined over an open dense subset
$U \subset X_0$ and such that
$f_0(\sigma(u)) = u$ for any $u \in U$.
By definition the points of $\cP_0$ are
 the lines $l$ that lie on the quadrics $Q_t, t \in \PP^3_0$. Denote by $l_u \in \cP_0$ the line $l_u = \sigma(u)$
for  points $u \in U$, i.e.
$$
\sigma\colon U \rightarrow \cP_0, \ x \mapsto l_u.
$$
Let $\pi\colon X \rightarrow \PP^3$ be the double covering,
and let  $i\colon X \rightarrow X$ be the involution
interchanging  two possibly coincident $\pi$-preimages of
the points $x \in \PP^3$.
Without any lost of generality  (e.g. by replacing $U$ by $U \cap i(U)$)
we may assume  that $U = i(U)$.
Let $D \subset W$ be  Zariski closure of  set
$$
\{ l_u\cap l_{i(u)}: u \in U \ \mbox{ and } \ u \not= i(u) \}.
$$
The variety $D$ is a 3-fold in $W$
that intersects the general quadric
$Q_x \subset \PP^3_x = x \times \PP^3$, $x = \pi(u)$
at {\it a unique} point --- the point $y(u)= l_u\cap l_{i(u)}$,
i.e. $DQ_x = 1$.

The 5-fold $W = (\PP^3 \times \PP^4) \cap H \cap (F(x;y) = 0)$
is an ample divisor in the 6-fold $Z = (\PP^3 \times \PP^3) \cap H$,
which in turn is an ample divisor in $\PP^3 \times \PP^4$.

Then by  Lefschetz hyperplane section Theorem
the restriction map defines  an isomorphisms
$$
H^4(\PP^3 \times \PP^4, \ZZ) \rightarrow H^4(Z, \ZZ) \rightarrow H^4(W, \ZZ).
$$
In particular, the codimension two subvariety $D \subset W$
is  a restriction  of a codimension two
subvariety of $\PP^3 \times \PP^4$ to $W$.

In the Chow ring
$$
A_*(\PP^3 \times \PP^4) = \ZZ[h_1,h_2]/(h_1^4,h_2^5),
$$
the class of the fibre $Q_x$ of
$p\colon W \rightarrow X$ is $2h_1^3h_2^2$. Since codimension 2 cycles on $\PP^3 \times \PP^4$
are generated over $\ZZ$ by $h_1^2, h_2^2$, and $h_1h_2$,
then the intersection number of any codimension 2 cycle on  $W$
with  general quadric $Q_x$ is {\it even}, which
contradicts  equality $DQ_x = 1$.
\end{proof}

Notice also that the varieties $X_0$ and $\widetilde{X}$ fulfill
conditions (\ref{h3-conditions}) from \ref{ai-3},
so $Br(X_0)$ and $Br(\widetilde{X})$ are isomorphic to
$H^3(X_0,\ZZ)$ and $H^3(\widetilde{X},\ZZ)$.

\begin{theorem}\label{non-zero-torsion}
The $\PP^1$-bundle $\cP_0$ represents a non-zero 2-torsion element
in  $Br(X) = H^3(\widetilde{X},\ZZ)$.
In particular, $\widetilde{X}$ and hence $X$ is non-rational.
\end{theorem}

\begin{proof}
Let $E_i, i = 1,\ldots,35$ be the exceptional divisors
of the blowup $\widetilde{X} \rightarrow X$ at the nodes $p_1,\ldots,p_{35}$.
Then by \cite{Gr}, for the Brauer groups of
$X_0 = X - \{ p_1,\ldots,p_{35} \} \cong \widetilde{X} - \cup \{ E_i : i = 1,\ldots, 35 \}$
there is an exact sequence
$$
0 \rightarrow Br(\widetilde{X}) \rightarrow Br(X_0)
   \rightarrow \bigoplus_{i=1}^{35} H^1(E_i,\QQ/\ZZ),
$$
and since for  surfaces $E_i \cong \PP^1 \times \PP^1$
one has $H^1(E_i,\QQ/\ZZ) = 0$, $i = 1,\ldots, 35$, then
$Br(\widetilde{X}) \cong Br(X_0)$.
\end{proof}

\medskip

It follows from Lemmas \ref{not-a-projectivization} and
\ref{ai-3} that $\cP_0$ represents
a non-zero 2-torsion element of $H^3(\widetilde{X},\ZZ)$.
Combining with  Lemma \ref{torsion-criterion}
we get   non-rationality of $\widetilde{X}$,
and hence --- {\it the non-rationality} of $X$.
This proves Proposition \ref{sds}. 

\section{Artin--Mumford quartic double solid}

\label{section:quartic}

\subsection{Quadrics in $\PP^3$ and Artin--Mumford quartic double solid}\label{ai-6}
Let $\PP^3 = \PP^3(y)$, $(y) = (y_0:\ldots:y_3)$
be the 3-dimensional complex projective space.
In the space $\PP^9 = \PP(H^0(\cO_{{\bf P}^3}(2))$ of quadrics
in $\PP^3$ regard the determinantals
$$
\Delta_1 \subset \Delta_2 \subset \Delta_3 \subset \PP^9
$$
where
$$
\Delta_k = \{ Q \in \PP^9: rank\ Q \le k, k = 1,2,3 \}.
$$
The elements of $\PP^9$ are  $\CC^*$-classes of
symmetric $4 \times 4$ matrices $Q = (q_{ij})$, $0 \le i,j \le 3$,
and the determinantals $\Delta_k$, $1 \le k \le 3$, defined by
vanishing $(k+1)\times (k+1)$ minors of $Q$ have the following
properties; for more details see e.g. \S 1 in \cite{Co}:

\medskip
{\it

$\Delta_3 \subset \PP^9$ is a quartic hypersurface;

\smallskip

$\Delta_2 = Sing\ \Delta_3$ has dimension 6 and degree 10;

\smallskip

$\Delta_1 = Sing\ \Delta_2 = v_2(\PP^3)$ is the Veronese image of $\PP^3$ in $\PP^9$;

\smallskip

The determinantal quartic $\Delta_3$
has an ordinary double singularity along $\Delta_2 - \Delta_1$.

}
\medskip

Consider   {\it general} $3$-space $\PP^3 = \PP^3(x) \subset \PP^9$.
As it follows from  previous considerations:

\medskip

{\it
$S = \PP^3 \cap \Delta_3$ is a quartic surface with only singularities
--- the 10 points of intersection
$\delta = \PP^3 \cap \Delta_2 = \{ p_1,\ldots,p_{10} \}$,
and any $p_k$, $k = 1,\ldots,10$, is an ordinary double point (a node) of $S$.\footnote{
The  quartic surfaces defined as  determinantal loci of 3-spaces
of quadrics in  projective 3-space appear in the works of A.\,Cayley
since the 80's of  19-th century under the name {\it quartic symmetroids}.
}
}

\medskip

Since $deg\,(S) = 4$ is an even number,  there exists a double
covering
$$
\pi\colon X \rightarrow \PP^3
$$
branched along $S$, i.e. $X$ is a determinantal quartic double solid.

The double solid $X$ has 10 nodes ---  isomorphic preimages of the 10
nodes $p_1,\ldots,p_{10}$ of  branch locus $S$, which we also denote
by $p_1,\ldots,p_{10}$. Let $\widetilde{X}$ be the blowup of $X$ at these 10 points. In the same way as in Section \ref{ai-4} we get:

\begin{proposition}
The group $H^3(\widetilde{X},\ZZ)$ contains a non-zero 2-torsion element;
in particular $X$ it is non-rational.
\end{proposition}

\begin{remark}
In  \cite{AM}, Artin and Mumford prove
stronger result: \ $Tors\,(H^3(X,\ZZ)) = \ZZ/2\ZZ$, by using
splitting of  discriminant curve for  natural conic bundle
structure on $X$, see also Theorem~2 in~\cite{Z}.
\end{remark}

\subsection{Artin--Mumford quartic double solids and Enriques surfaces}\label{ai-8}

We start by recalling   well known connection
between Artin--Mumford double solids and Enriques surfaces, defined by Reye
congruences, see e.g. \cite{Co}.
In  above notation,
the Artin--Mumford double solids are defined by the general 3-spaces
$\PP^3(x)$ in the space
$\PP^9 = \PP(H^0(\cO_{\PP^3(y)}(2))$
of quadrics in $\PP^3(y)$, $(y) = (y_0:\ldots:y_3)$.
Let
$$
\{ Q_x \} = \{Q_x \subset \PP^3(y): x \in \PP^3(x) = \PP^3(x_0:\ldots:x_3) \}
$$
be the set of quadrics in $\PP^3(y)$ defined by the 3-space $\PP^3(x)$.
Let $G$ be the Grassmannian of lines $l \subset \PP^3(y)$.

It is known that  general line $l \subset \PP^3(y)$ lies
on a unique quadric from the family $\{ Q_x \}$, and the set of lines
$$
R = \{ l \in G : \mbox{ the line } l \subset \PP^3(y) \mbox{ lies in a }
\PP^1\mbox{-family of quadrics } Q_x \}
$$
is an Enriques surface in $G = G(2,4)$
called classically a {\it Reye congruence}, see \cite{Co}.
Let ${\tau}$ be an involution
$$
(x,y) \stackrel{{\tau}}\longleftrightarrow (y,x)
$$
on $\PP^3(x) \times \PP^3(y)$.
The fixed point set of ${\tau}$ is the diagonal
$\Delta$ defined by $\{x = y\}$ in $\PP^3(x) \times \PP^3(y)$.
For a quadratic form
$$
Q(y) = \sum_{0 \le i,j \le 3} q_{ij}y_iy_j, \ q_{ij} = q_{ji},
$$
let
$$
B(x,y) = \sum_{0 \le i,j \le 3} q_{ij}x_iy_j
$$
be its corresponding bilinear form. Then a basis $Q_0(y), \ldots , Q_3(y)$
of $\PP^3(x) \subset \PP^9$ defines a quadruple of bilinear forms
$B_0(x,y), \ldots, B_3(x,y)$, and hence --- a linear section
$$
\widetilde{S} = (\PP^3(x) \times \PP^3(y))\cap H_0 \cap \ldots \cap H_3
$$
where $H_i = (B_i(x,y) = 0)$.
For a  general choice of
$$\PP^3(x) = \langle Q_0,\ldots,Q_3\rangle$$
the set $\widetilde{S}$ is a smooth complete intersection of 4
hyperplane sections of $\PP^3(x) \times \PP^3(y)$, and hence
$\widetilde{S}$ is a smooth K3 surface --- {\it  Steiner} K3 surface
 in  3-space of quadrics $\PP^3(y)$.
Since all $B_i$ are invariant under the involution ${\tau}$, then
$\widetilde{S}$ is also invariant under ${\tau}$, i.e.
${\tau}(\widetilde{S}) = \widetilde{S}$.
Therefore ${\tau}$ restricts to an involution
${\tau}\colon \widetilde{S} \rightarrow \widetilde{S}$;
and since for  general $\PP^3(x)$ the
surface $\widetilde{S}$ does not intersect  diagonal $\Delta$
we conclude  ${\tau}$ is without fixed points on $\widetilde{S}$.
The K3 surface $\widetilde{S}$ has  following properties
(see \cite{Co}, \cite{O}):

\medskip

Let $\PP^3(x)$ be a general 3-space of the 9-space $\PP^9$ of quadrics in
$\PP^3(y)$, and let $S = D_3 \subset \PP^3(x)$, $R \subset G(2,4)$ and
$\widetilde{S}$ be correspondingly the quartic symmetroid,
the Enriques surface (the Reye congruence),
and the Steiner K3 surface defined by $\PP^3(x)$. Then:

\medskip

(i)
{\it $\widetilde{S}$ is the blowup of $S$
at its 10 nodes $\delta = \{ p_1,\ldots,p_{10} \}$};

\medskip

(ii)
{\it $R \subset G = G(2,4)$ is isomorphic to the quotient
$\widetilde{S}/\tau$ of $\widetilde{S}$ by the involution $\tau$}.

\medskip

Let $\pi\colon X \rightarrow \PP^3(x)$ be the Artin--Mumford double solid,
defined by the general 3-space $\PP^3(x) \subset \PP^9$,
let $G = G(1:{\PP^3}(y))$ be as above,  and let
$$
\widetilde{G} = \{ (x,l) \in \PP^3(x) \times G : l \subset Q_x \}.
$$
Then (see \S 9 in \cite{Be}):

\medskip

(iii)
{\it $\widetilde{G} = \cP$ (see the proof of Proposition \ref{sds}),
and the projection $\widetilde{G} \rightarrow G$, $(x,l) \mapsto l$
is a blowup of the Enriques surface $R \subset G = G(2,4)$}.

\medskip

(iv)
{\it The projection $\sigma\colon \widetilde{G} \rightarrow \PP^3$, $(x,l) \mapsto x$
factorizes into
$$
\widetilde{G} \stackrel{f}\longrightarrow X \stackrel{\pi}\longrightarrow \PP^3(x),
$$
and the restriction $\widetilde{G}_0 \rightarrow X_0$ of $f$ over $X_0 \subset X$
coincides with the $\PP^1$-bundle $f_0\colon \cP_0 \rightarrow X_0$:
}

%
%
%
%


$$
\xymatrix{
\cP_0\ \subset \ \cP\ \cong \ \widetilde{G}\ar@<-6ex>[d]_{f_0}\ar@<0ex>[d]^{f}\ar[r]^-\sigma & G(2,4)\supset R\\
\!\!\!\!\!\!\!\!\!\!\!\!\!\!\!\!X_0\ \subset \ X \ar@<-6ex>[d]_{\pi_0}\ar@<0ex>[d]^{\pi}& \\
\!\!\!\!\!\!\!\!\!\!\!\!\!\!\!\!\PP^3\ \subset \ \PP^3\\
}
$$

\subsection{The non-rationality of $X$ by the Criterion \ref{torsion-criterion} (see \cite{Be})}

We observe  that since   $\sigma\colon \cP  = \widetilde{G} \rightarrow G(2,4)$
is a blowup of the surface $R$ in the 4-fold $G(2,4)$,
then
$$H^4(\cP, \ZZ) =  \sigma^*H^4(G(2,4),\ZZ) \oplus \sigma^{-1}H^2(R,\ZZ)$$
$$\cong H^4(G(2,4),\ZZ) \oplus H^2(R,\ZZ).$$

Furthermore since $R$ is an Enriques surface, then
$c_1(R) \in H^2(R,\ZZ)$ is an element
of order 2. Therefore $Z = \sigma^{-1}c_1(R)$ is an element
of order 2 in $H^4(\cP, \ZZ)$.
After restriction, we get an element $Z_0 \in H^4(\cP_0,\ZZ)$ of order 2.

Since $f_0\colon\cP_0 \rightarrow X_0$
is a $\PP^1$-bundle, then all  fibers of
$f_0$ are isomorphic to 2-dimensional spheres
$S^2$. Therefore the integral  cohomology of
$\cP_0$ and $X_0$ fit in the  Gysin sequence
for $S^2$-fibration:
$$
\ldots \longrightarrow H^3(\cP_0,\ZZ) \longrightarrow H^1(X_0,\ZZ)
\stackrel{e}\longrightarrow H^4(X_0,\ZZ)
$$
$$
\stackrel{f_0^*}\longrightarrow H^4(\cP_0,\ZZ)
\stackrel{f_{0*}}\longrightarrow H^2(X_0,\ZZ) \longrightarrow \ldots.
$$

Here  $e$ is the cup-product with the Euler class $e(f_0) \in H^3(X_0,\ZZ)$
of $f_0$, see Chapter III, \S 14 in \cite{BT} and 4.11 in \cite{Hi}.

If $Im(e) \not=0$ then any non-zero element of $Im(e) \in H^4(X_0,\ZZ)$ is a 2-torsion
element, since $2e = 0$, see Theorem 4.11.2 (I) in \cite{Hi}.

In  case when $Im(e) =0$ then $Z_0 \in Tors_2\,(H^4(\cP_0,\ZZ))$ must be an image
$Z_0 = f_0^*(C_0)$ of an element $C_0 \in H^4(X_0,\ZZ)$,
since $Tors\,(H^2(X_0,\ZZ))$ = $0$ (see p. 30 in \cite{Be}). Since in this case $f_0^*$ is an embedding and $Z_0$ is a non-zero 2-torsion
element of $H^4(\cP_0,\ZZ)$,
then $C_0$ is also a non-zero 2-torsion element of $H^4(X_0,\ZZ)$.

Thus in both cases there exists a 2-torsion element $C_0 \in  H^4(X_0,\ZZ)$.

Let $\sigma_X\colon \widetilde{X} \rightarrow X$ be the blowup of $X$ at
the 10 nodes $p_1,\ldots,p_{10}$ of $X$, and let
$E_i = \sigma_X^{-1}(p_i) \cong \PP^1 \times \PP^1$
be the 10 exceptional divisors on $\widetilde{X}$.
Since $\widetilde{X}$ is isomorphic to a disjoint union of $X_0$
and $E_i$, $i = 1,\ldots,10$, and
$H^3(E_i,\ZZ) = H^3(\PP^1 \times \PP^1,\ZZ) = 0$,
then $H^4(X_0,\ZZ)$ is embedded isomorphically in
$H^4(\widetilde{X},\ZZ)$. In particular
$C_0 \in  H^4(X_0,\ZZ)$ is embedded as an element
$C$ of order two in $H^4(\widetilde{X},\ZZ)$.

Since for a smooth projective complex
threefold $\widetilde{X}$ one has
$$Tors\,(H^4(\widetilde{X},\ZZ)) \cong Tors\,(H^3(\widetilde{X},\ZZ))$$
(see \S 1 in \cite{AM}), the 2-torsion element $C \in H^4(\widetilde{X},\ZZ)$
represents equally a 2-torsion element $Z_C \in H^3(\widetilde{X},\ZZ)$.
By the Criterion \ref{torsion-criterion} the last yields that
$\widetilde{X}$ (and hence $X$) is non-rational.

\begin{remark}
\label{remark:quadric solid}
It was suggested to us by  K.\,Shramov that  methods of~\cite{AMG} can be applied to  double covering of quadric ramified in octic with 20 singular points. More precisely we  consider
a divisor of bidegree (1,2) in $Q\times \PP^3$, where $Q$ is a quadric threefold. In this case we get a quadric fibration
given by a map $\cO_Q(-1)\to S^2(E^*)$, where $E$ is a trivial vector bundle of rank 4.
We get a 2-torsion (and hence nonrationality) in a middle cohomology of a double quadric with 20 nodal singular points. Using the fact that  double covering of quadric ramified in octic with 20 singular points is a degeneration of three dimensional quartic we will study its Landau--Ginzburg model in Section 5.

\end{remark}

\section{Mirror Side}

In this section we turn to Homological Mirror Symmetry in an attempt  to show that phenomena  observed  in previous sections is a part of much more general scheme.
We briefly outline in  Figure \ref{ClassicalHMS}  a schematic picture of classical Homological
Mirror Symmetry, in a version relevant for our purpose.
For more details see \cite{KRIS1}.

\begin{figure}
\begin{center}
\begin{tabular}[c]{c|c} \hline
A-models (symplectic) & B-models (algebraic)
\\ \hline \\
$ X = (X,\omega) $ a closed symplectic manifold &
$ X $ a smooth projective variety
\\[1.2em]
\begin{minipage}[t]{0.45\linewidth}
\emph{Fukaya category} $ \Fuk(X)$.
Objects are Lagrangian submanifolds $L$ which may be equipped with flat line
bundles.
Morphisms are given by Floer cohomology $ HF^*(L_0,L_1) $.
\end{minipage} &
\begin{minipage}[t]{0.45\linewidth}
\emph{Derived category} $ D^b(X) $.
Objects are complexes of coherent sheaves $ \mathcal{E} $.
Morphisms are $ Ext^*({\mathcal E}_0,{\mathcal E}_1) $.
\end{minipage}
\\
\multicolumn{2}{c}{$$\xymatrix{ \ar@{<->}[drr]&&\ar@{<->}[dll]\\ && }$$}
\medskip
\\
\begin{minipage}[t]{0.45\linewidth}
A non-compact symplectic manifold $Y$ with a proper map
$ W\colon Y \to \CC $ which is a symplectic fibration with singularities.
\end{minipage}\ & \
\begin{minipage}[t]{0.45\linewidth}
$ Y $ a smooth quasi-projective variety with a proper holomorphic map
$ W\colon Y \to \CC $.
\end{minipage}
\\[3.5em]
\\
\begin{minipage}[t]{0.45\linewidth}
\emph{Fukaya--Seidel category of the Landau--Ginzburg model} $ FS(LG(Y))$:
Objects are Lagrangian submanifolds $ L \subset Y $ which, at infinity, are
fibered over $ \RR^+ \subset \CC $.
The morphisms are $ HF^*(L_0^+,L_1) $, where the superscript
$ + $ indicates a perturbation removing intersection points at infinity.
\end{minipage} &
\begin{minipage}[t]{0.45\linewidth}
The category $ D^b_{sing}(W) $ of algebraic $ B $-branes which is obtained by
considering the singular fibers $ Y_z = W^{-1}(z) $, dividing
$ D^b(Y_z) $ by the subcategory of perfect complexes {\it Perf}$ \,(Y_z) $,
and then taking the direct sum over all such $ z $.
\end{minipage}
\end{tabular}
\end{center}
\caption{Classical Homological Mirror Symmetry.}
\label{ClassicalHMS}
\end{figure}


In what follows we describe 
fiberwise compactifications of weak Landau--Ginzburg models of quartic double solid, Fano threefold $V_{10}$, and of sextic double solid  (see~\cite{Prz09a} and~\cite{Prz13}).
We conjecture that these compactifications are Landau--Ginzburg models of the Artin--Mumford example, $V_{10}$, and sextic double solid
correspondingly in the sense of HMS.

Throughout this section we use the following standard notations for blowup.
Consider  affine variety
$$
\{F(x_1,\ldots,x_n)=0\}\subset \Aff(x_1,\ldots,x_n).
$$

We blow up  affine space $\{x_1=\ldots=x_k=0\}$.
The blown up hypersurface is given by the system of equations
$$
\left\{%
\begin{array}{l}
    F(x_1,\ldots,x_n)=0,\\
    x_ix_j^\prime=x_i^\prime x_j, \ \ 1\leq i,j\leq k,
\end{array}%
\right.
$$
in
$$
\Aff(x_1,\ldots,x_n)\times \PP(x^\prime_1:\ldots:x_n^\prime).
$$
Consider  local chart  $x_1^\prime\neq 0$.
We choose coordinates
$$
x_1,\frac{x_2^\prime}{x_1^\prime},\ldots,\frac{x_k^\prime}{x_1^\prime},x_{k+1},\ldots,x_n.
$$
In these coordinates blown up variety is zero locus of  polynomial given by division of
$$
F(x_1,x_1x_2^\prime,\ldots,x_1x_k^\prime,x_{k+1},\ldots,x_n)
$$
by maximal possible power of $x_1$. We use notations $x_i$'s for coordinates in this local chart instead of $\frac{x_i^\prime}{x_1^\prime}$'s
for simplicity. We denote this local chart by $x_1\neq 0$.

We embed  fiberwise above  pencil in a projective space  or product of projective spaces  and  then resolve singularities. 
All Calabi--Yau compactifications (see~\cite{Prz09a}) are birational in codimension one.

\subsection{The Landau--Ginzburg model of quartic double solid.}
\label{LG for quartic2solid}

The weak Landau--Ginzburg model for quartic double solid is given by

$$
f=\frac{(x+y+1)^4}{xyz}+z\ \in \CC [x^{\pm 1}, y^{\pm 1}, z^{\pm 1}].
$$

We compactify  pencil $\{f=\lambda$, $\lambda\in \CC\}$, in the
neighborhood of $\lambda=0$ in $\PP(x:y:z:t)\times \Aff(\lambda)$ and get  hypersurface
$$
\{(x+y+t)^4+xyz(z-\lambda t)=0\} \ \subset \PP(x:y:z:t)\times
\Aff(\lambda).
$$
Its singularities are seven lines
$$
l_0=\{x+y+t=z=\lambda=0\},\ \ \
l_1=\{x=y=t=0\},\ \ \
l_2=\{x+y=z=t=0\},
$$
$$
l_3=\{x=y+t=z=0\},\ \ \
l_4=\{x=y+t=z+\lambda y=0\},\ \ \
l_5=\{x+t=y=z=0\},
$$
$$
l_6=\{x+t=y=z+\lambda x=0\}.
$$

\begin{figure}[htbp]
  \begin{center}
\includegraphics[width=0.3\textwidth]{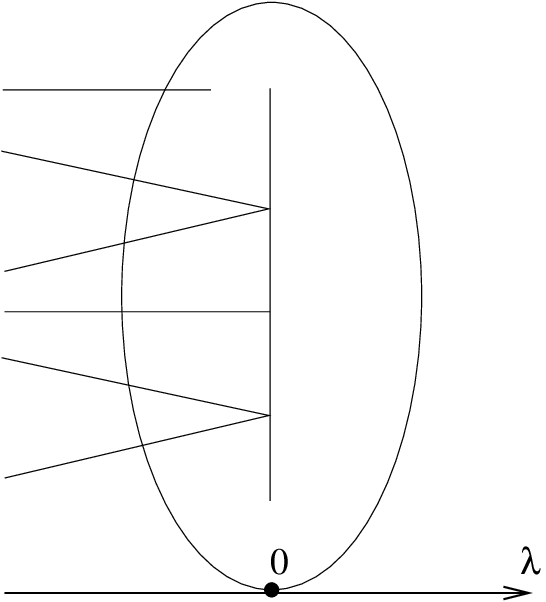}
\caption{Singularities for quartic double
solid.}\label{figure:quart_sol_1}
  \end{center}
\end{figure}

Generically above  singularities  are locally products
of du Val singularities of type $A_3$ by affine line. ``Horizontal''
lines $l_2$--$l_6$ intersect ``vertical'' line $l_0$; moreover,
pairs of lines $l_3$ and $l_4$, $l_5$ and $l_6$ intersect $l_0$ at
one point (see Figure~\ref{figure:quart_sol_1}).

We resolve singularities by blowing up these lines.  At first we blow up the vertical line $l_0$ twice.
After this the singularities are proper transforms of
lines $l_1$--$l_6$ and five lines lying on the exceptional divisors. Each of them intersect proper transform of one of lines
$l_2$--$l_6$. After blowing up these five lines we get threefold with six
lines of singularities coming from $l_1$--$l_6$ which are of type
$A_3$ along a horizontal affine line globally. Blowing them up
fiberwise we get the final resolution. We carry this procedure in the following steps:

\medskip

{\bf Step 0.} The line $l_1$ is of type $A_3$ along affine line
globally. Blowing it up twice we get horizontal exceptional
fibers, so they do not give an additional component for fiber over
$\lambda=0$. We proceed resolution in the neighbourhood of line
$l_0$.

\medskip

{\bf Step 1.} Let $a=x+y+t$. Then our variety is given by
$$
\{a^4+xyz^2=\lambda xyz (a-x-y)\}\ \subset \PP(x:y:z:a)\times
\Aff(\lambda)
$$
and $l_0=\{a=z=\lambda=0\}$. There are two similar local charts: $x\neq
0$ and $y\neq 0$. Consider  local chart $y\neq 0$. It contains lines of singularities
$l_0$, $l_2$--$l_4$. We study the resolution in this chart
and double the picture over lines $l_3$, $l_4$. In this local chart we
have an affine hypersurface
$$
a^4+xz^2=\lambda xz (a-x-1)
$$
and we need to blow up  line $l_0=\{a=z=\lambda=0\}$.

{\bf The local chart 1a: $a\neq 0$.} We have hypersurface
$$
a^2+xz^2=\lambda xz (a-x-1).
$$
The exceptional divisor is given by equation $a=0$, so it consists
of three components
$$
E_1^a=\{a=x=0\}, \ \ \ E_2^a=\{a=z+(x+1)\lambda=0\}, \ \ \
E_3^a=\{a=z=0\}.
$$
The proper transform of the fiber over $\lambda=0$ is
$E_0=\{\lambda=a^2+xz^2=0\}$. The singularities are:
$$
l_1^a=\{x=z=a=0\}, \ \ l_2^a=\{x=\lambda+z=a=0\},
$$
$$
l_3^a=\{z=a=\lambda=0\}, \ \ l_4^a=\{x+1=z=a=0\}.
$$
We have:
$$
E_2^a\cap E_3^a=l_3^a\cup l_4^a,\ \ \  E_1^a\cap E_3^a=l_1^a, \ \ \
E_0\cap E_2^a\cap E_3^a=l_3^a.
$$
All proper transforms of lines $l_2$--$l_6$ do not lie in this chart.

\medskip

{\bf The local chart 1z: $z\neq 0$.} There is nothing new in this chart:
all we are interested in is contained in the chart {\bf 1a}.

\medskip

{\bf The local chart 1$\boldsymbol{\lambda}$: $\lambda\neq 0$.} We have
hypersurface
$$
\lambda^2a^4+xz^2= xz (\lambda a-x-1).
$$
The exceptional divisor is given by equation $\lambda=0$, so it
consists of three components
$$
E_1^\lambda=\{\lambda=x=0\},  \ \ \ E_2^\lambda=\{\lambda=z+x+1=0\},\ \
\ E_3^\lambda=\{\lambda=z=0\}.
$$
The proper transform of  fiber over $\lambda=0$ does not lie in
this chart. We have:
$$
E_1^\lambda=E_1^a,\ \ \ E_2^\lambda=E_2^a, \ \ \ E_3^\lambda=E_3^a.
$$
The singularities are:
$$
l_1^\lambda=\{x=z=\lambda=0\}=E_1^\lambda\cap E_3^\lambda,
$$
$$
l_2^\lambda=\{a=x=z=0\} \ \ \  \mbox{--- proper transform of $l_3$},
$$
$$
l_3^\lambda=\{x+1=z=\lambda=0\}=E_2^\lambda \cap E_3^\lambda,
$$
$$
l_4^\lambda=\{a=z=x+1=0\} \ \ \  \mbox{--- proper transform of
$l_2$},
$$
$$
l_5=\{x=z+1=\lambda=0\}=E_1^\lambda\cap E_2^\lambda,
$$
$$
l_6^\lambda=\{x=z+1=a=0\} \ \ \  \mbox{--- proper transform of
$l_4$}.
$$

So, after  first blow-up we get a configuration of components of  central fiber drawn on Figure~\ref{figure:first blowup}.

\begin{figure}[htbp]
  \begin{center}
\includegraphics[width=0.3\textwidth]{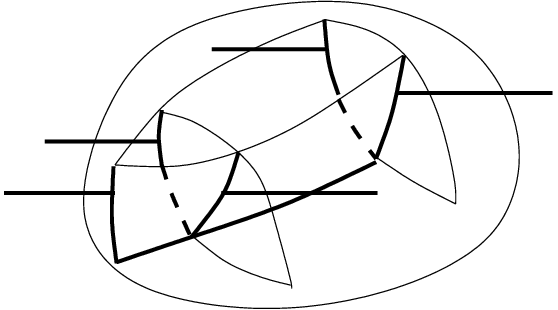}
\caption{The picture after the first blowup.}\label{figure:first
blowup}
  \end{center}
\end{figure}

Then we blow up the line $l_3^a$. It is enough to consider it in the chart
{\bf 1a}. That is, we blow up the line
$$
\{z=a=\lambda=0\}
$$
at
$$
\{a^2+xz^2-\lambda xz(a-x-1)=0\}.
$$

%
%
%
%
%
%
%

The only meaningful  local chart is $\lambda\neq 0$.
In this chart get the
hypersurface
$$
\{a^2+xz^2-xz(\lambda a-x-1)=0\}.
$$
The exceptional divisor is
$$
E^{a,\lambda}=\{\lambda=a^2-xz(z+x+1)=0\}.
$$
The singularities in its neighborhood are
$$
\{a=x=z=0\} \ \ \  \mbox{--- proper transform of $l_1^a$},
$$
$$
\{a=x=z+1=0\} \ \ \  \mbox{--- proper transform of $l_2^a$},
$$
$$
\{a=x+1=z=0\} \ \ \  \mbox{--- proper transform of $l_4^a$}.
$$
All of them lie on the exceptional divisor. So we did not get
``new'' singularities after this blowup. The divisors $E_2^a$,
$E_3^a$ now intersect only by proper transform of $l_4^a$; the
divisor $E_1^a$ intersect $E_2^a$ and $E_3^a$ in two separated lines
both intersect the line $E_1^a\cap E^{a,\lambda}=\{x=a=\lambda=0\}$.
The proper image of $E_0$ intersects only $E^{a,\lambda}$ by a line
lying far from the rest exceptional divisors.

\medskip

Now we  blow up  line $l_4^a=l_3^\lambda$. The line $l_3^a$ does
not lie at the chart {\bf 1$\boldsymbol{\lambda}$} so we can consider
this blowup only in the chart {\bf 1$\boldsymbol{\lambda}$}. We make a change
 of
variables $x\to x-1$. Then we get a  hypersurface
$$
\{\lambda^2 a^4+(x-1)z^2=(x-1)z(\lambda a -x)=0\}
$$
and then we need to blow up the line
$$
\{x=z=\lambda=0\}.
$$

We get one exceptional divisor, proper images of lines $l_5^\lambda$ and $l_6^\lambda$ that lie far from exceptional divisor,
proper images of $l_1^\lambda$ and $l_2^\lambda$ (we will discuss them later) and a proper image of $l_4^\lambda$ (in the other words, of $l_2$).
It is globally of type $A_3$ along a line, so it resolves
horizontally and does not give an exceptional divisors over
$\lambda=0$.

So, after this blowup the divisors $E_2^\lambda$ and $E_3^\lambda$
are separated.

\medskip

Now we  blow up the line $l_1^\lambda$. As before, we can do it in the
chart {\bf 1$\boldsymbol{\lambda}$.} We have a hypersurface
$$
\{\lambda^2 a^4+xz(z+x+1-\lambda a)=0\}
$$
and need to blow up the line $\{x=z=\lambda=0\}$.

%
%
%
%
%
%

We get proper transforms of $l_3^\lambda$ and $l_4^\lambda$ we already discussed,
proper transforms of $l_5^\lambda$ and $l_6^\lambda$ we mention in the next paragraph,
and a proper transform of $l_2^\lambda$ (in the other words, of $l_3$). It is
globally of type $A_3$ along a line, so it resolves horizontally
and does not give an exceptional divisors in the central fiber.

\medskip

Finally, the picture we get after blowing up the line $l_5^\lambda$ is
very similar to the picture we get after blowing up the line
$l_1^\lambda$.

\begin{figure}[htbp]
\begin{center}
\includegraphics[width=.5\textwidth]{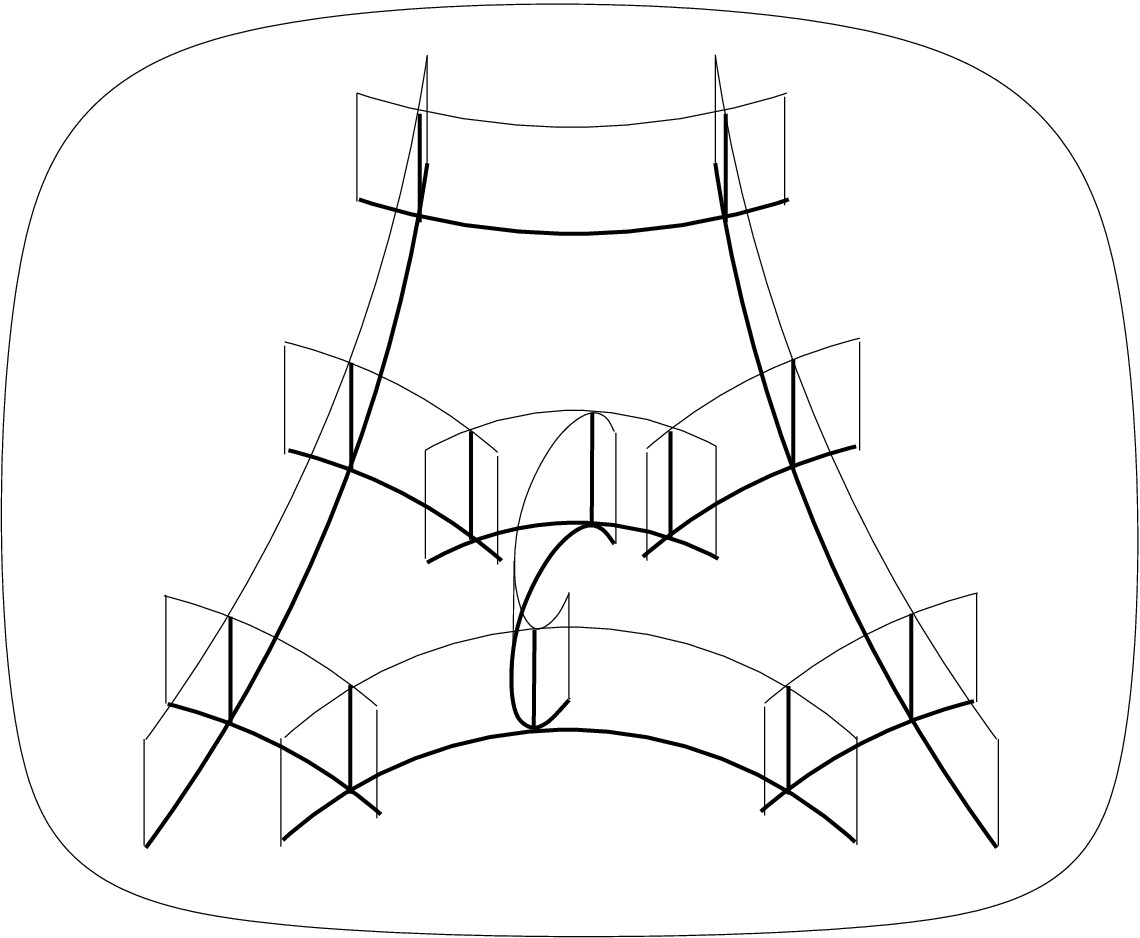}
\end{center}
\caption{The fiber over $0$ in Landau--Ginzburg model for a quartic
double solid.} \label{figure:quartic double solid}
\end{figure}

We summarize  final picture of resolved singularities (see
Figure~\ref{figure:quartic double solid}).

\medskip

Via direct calculations (see \cite{MARNON}, \cite{KNS}) we get

\begin{proposition}  The monodromy of the singular fiber at zero of the  Landau--Ginzburg model for a quartic
double solid with 10 singular points is strictly unipotent.
\end{proposition}

The proof of above proposition is based on  analysis of monodromy  change
under conifold transition.

\subsection{The Landau--Ginzburg model of $V_{10}$.}
\label{LG for V10}

The weak Landau--Ginzburg model for a Fano variety $V_{10}$ is
$$
f=\frac{(x^2+x+y+z+xy+xz+yz)^2}{xyz}\ \in \CC [x^{\pm 1}, y^{\pm 1}, z^{\pm 1}].
$$

\begin{figure}[htbp]
  \begin{center}
\includegraphics[width=0.4\textwidth]{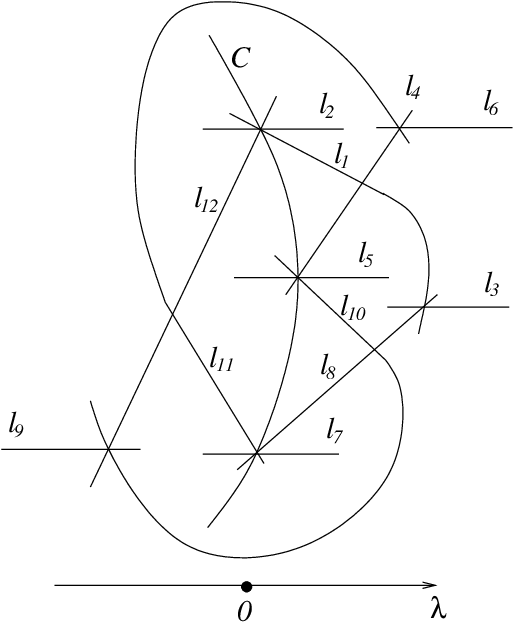}
\caption{Singularities for $V_{10}$.}\label{figure:V10_sing}
  \end{center}
\end{figure}

Compactifying  the pencil $\{f=\lambda$, $\lambda\in \CC\}$, in the
neighborhood of $\lambda=0$ in $\PP(x:y:z:t)\times \Aff(\lambda)$
we get a hypersurface
$$
\{(x^2+xt+zt+xz+yt+yz+xy)^2=\lambda xyzt
\} \ \subset \PP(x:y:z:t)\times
\Aff(\lambda).
$$
Its singularities are twelve lines
$$
l_1=\{x+z=t=\lambda=0\},\ \ \
l_2=\{x=z=t=0\},\ \ \
l_3=\{x+z=y=t=0\},
$$
$$
l_4=\{x+y=t=\lambda=0\},\ \ \
l_5=\{x=y=t=0\},\ \ \
l_6=\{x+y=z=t=0\},
$$
$$
l_7=\{x=y=z=0\},\ \ \
l_8=\{x+z=y=\lambda=0\},\ \ \
l_9=\{x+t=y=z=0\},
$$
$$
l_{10}=\{x=y+t=z=0\},\ \ \
l_{11}=\{x+t=y=\lambda=0\},\ \ \
l_{12}=\{x+t=z=\lambda=0\},
$$
and a conic
$$
C=\{x=yt+zt+yz=\lambda=0\}
$$
(see Figure~\ref{figure:V10_sing}).

There is a symmetry $x \leftrightarrow y \leftrightarrow z$, so we have three types of singular lines: two horizontal line types and
one vertical line type.

\begin{figure}[htbp]
  \begin{center}
\includegraphics[width=0.3\textwidth]{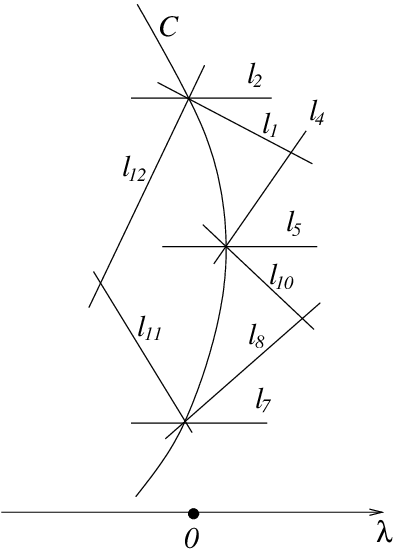}
\caption{Singularities for $V_{10}$.}\label{figure:V10_first_blowup}
  \end{center}
\end{figure}

We blow up $l_6$, we put $a=x+y$, and consider a local chart $x=1$.
In this chart  coordinates of  our family can be written as
$$
\{(a+at+zt+az)^2=\lambda (a-1)zt.
\}
$$
In the neighborhood of $l_6$ it is analytically equivalent to a hypersurface
$\{a^2=\lambda zt\}$.
In this local chart $l_6$, $l_{11}$, and $l_4$ are given by equations
$a=z=t=0$, $a=z=\lambda=0$, and $a=t=\lambda=0$ correspondingly.
They are intersecting transversally lines of singularities of type $A_1$.
So, blowing $l_6$ up we get one horizontal exceptional divisor. In its neighborhood the singularities (proper images of $l_{11}$ and $l_4$)
are lines of singularities of type $A_1$. Similarly, by symmetry,  the same holds in a neighborhood of lines $l_3$ and $l_9$.
After blowups performed above the singularities can be seen  on Figure~\ref{figure:V10_first_blowup}.

Let us blow up $l_7$  in  local chart $t=1$.
We have a hypersurface
$$
\{(x^2+x+z+xz+y+yz+xy)^2=\lambda xyz
\}.
$$

Analytically,  in a neighborhood of $l_7$ it is isomorphic  to a hypersurface
$$
\{(x+y+z+yz)^2=\lambda xyz\}.
$$
The lines $l_7$, $l_{8}$, $l_{11}$ and $C$ are given by equations
$x=y=z=0$, $x+z=y=\lambda=0$, $x+y=z=\lambda=0$, and $y+z+yz=x=\lambda=0$ correspondingly.

Consider a local chart $x\neq 0$ in the above blowup.
We get a hypersurface
$$
\{(1+y+z+xyz)^2=\lambda xyz
\}.
$$
The exceptional divisor is given by
$$
\{x=y+z+1=0\}
$$
and singularities in its neighborhood are given by
$$
l^x_1=\{x=y=z+1=0\},\ \ \
l^x_2=\{x=y+1=z=0\},\ \ \
l^x=\{x=y+z+1=\lambda=0\},
$$
$$
l_8=\{y=z+1=\lambda=0\},\ \ \
l_{11}=\{y+1=z=\lambda=0\}.
$$

In  neighborhood of  $l^x_1$ the singularities are three intersecting lines:
one horizontal  $l^x_1$ and two  vertical lines $l^x$ and $l_8$. They are analytically equivalent to
singular lines on the hypersurface $\{a^2=\lambda xy\}$.
We blow op  $l^x_1$ first and then  $l^x$ and $l_8$. We get two non-intersecting exceptional divisors in the central fiber coming from
$l^x$ and $l_8$.


Consider now a local chart $y\neq 0$ in the blowup.
We get a hypersurface
$$
\{(x+1+z+yz)^2=\lambda xyz
\}.
$$
The exceptional divisor is given by
$$
\{y=x+z+1=0\}
$$
and singularities in its neighborhood are given by
$$
l^x=\{y=x+z+1=\lambda=0\},\ \ \
l^y=\{x=y=z+1=0\},\ \ \
l^x_2=\{y=z=x+1=0\},
$$
$$
l_{11}=\{x+1=z=\lambda=0\},\ \ \
C=\{x=1+z+yz=\lambda=0\}.
$$

In  neighborhood of $l^y$ the singularities form three intersecting lines of ordinary double points $l^y$, $l_{11}$, and $C$ as before so we can
resolve them in similar  way.

Finally, we repeat with no change procedure  in the last local chart $z\neq 0$.
The lines $l_1$ and $l_4$ intersects transversally and are of type $A_1$. Blowing them up one-by-one, we get in the central fiber two  exceptional divisors
 intersecting  in a line.

\begin{figure}[htbp]
  \begin{center}
\includegraphics[width=0.4\textwidth]{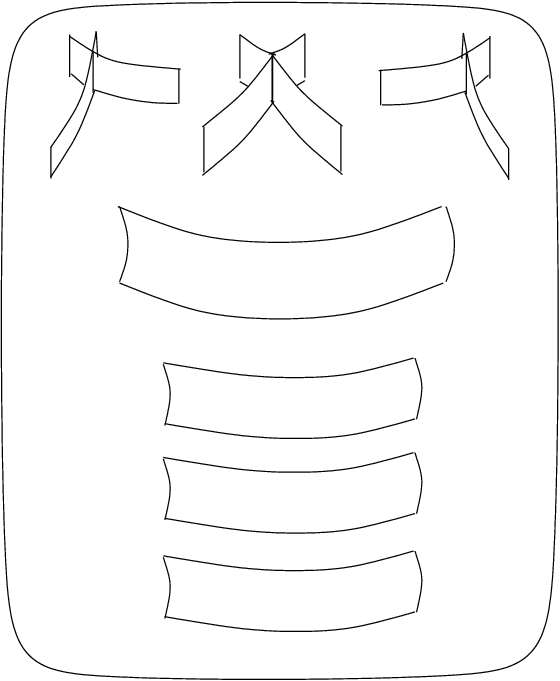}
\caption{Singularities for $V_{10}$.}\label{figure:V10_final}
  \end{center}
\end{figure}

The central fiber of  resolution shown  on Figure~\ref{figure:V10_final}.
There are eleven surfaces.

As before  direct calculations  based on  \cite{MARNON}, \cite{KNS} give:

\begin{proposition}  The monodromy of the singular fiber at zero of the  Landau--Ginzburg model for $V_{10}$ with  10 singular points is strictly unipotent.
\end{proposition}

%
%
%

\subsection{The Landau--Ginzburg model of sextic double solid.}

The weak Landau--Ginzburg model for a sextic double solid is

$$
\frac{(x+y+z+1)^6}{xyz}\in \CC [x^{\pm 1}, y^{\pm 1}, z^{\pm 1}].
$$

\begin{figure}[htbp]
\begin{center}
\includegraphics[width=.4\textwidth]{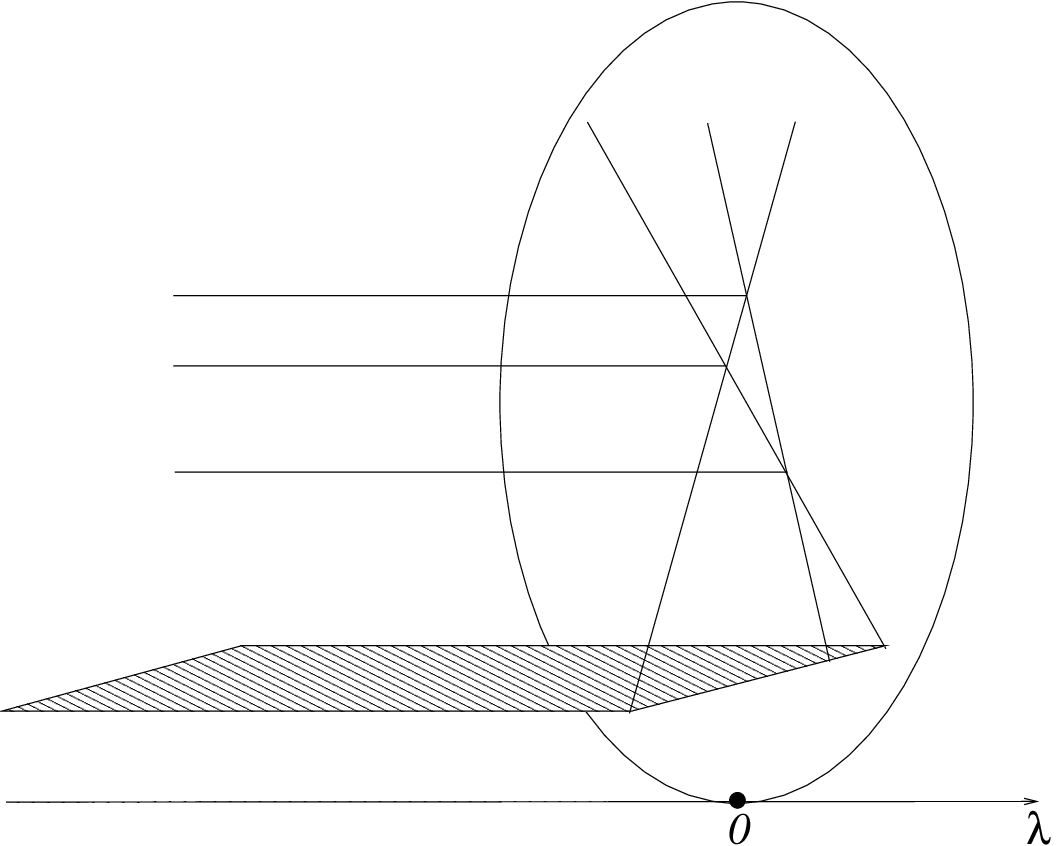}
\end{center}
\caption{The singularities for sextic double solid.}
\label{figure:sextic solid singularities}
\end{figure}

We are compactifying  it in a projective space.
The singularities are drawn on Figure~\ref{figure:sextic solid singularities}.
They are three vertical lines, three horizontal lines and a horizontal plane (lines are symmetric with respect to changing coordinates
$x\leftrightarrow y\leftrightarrow z$).

\begin{figure}[htbp]
\begin{center}
\includegraphics[width=.35\textwidth]{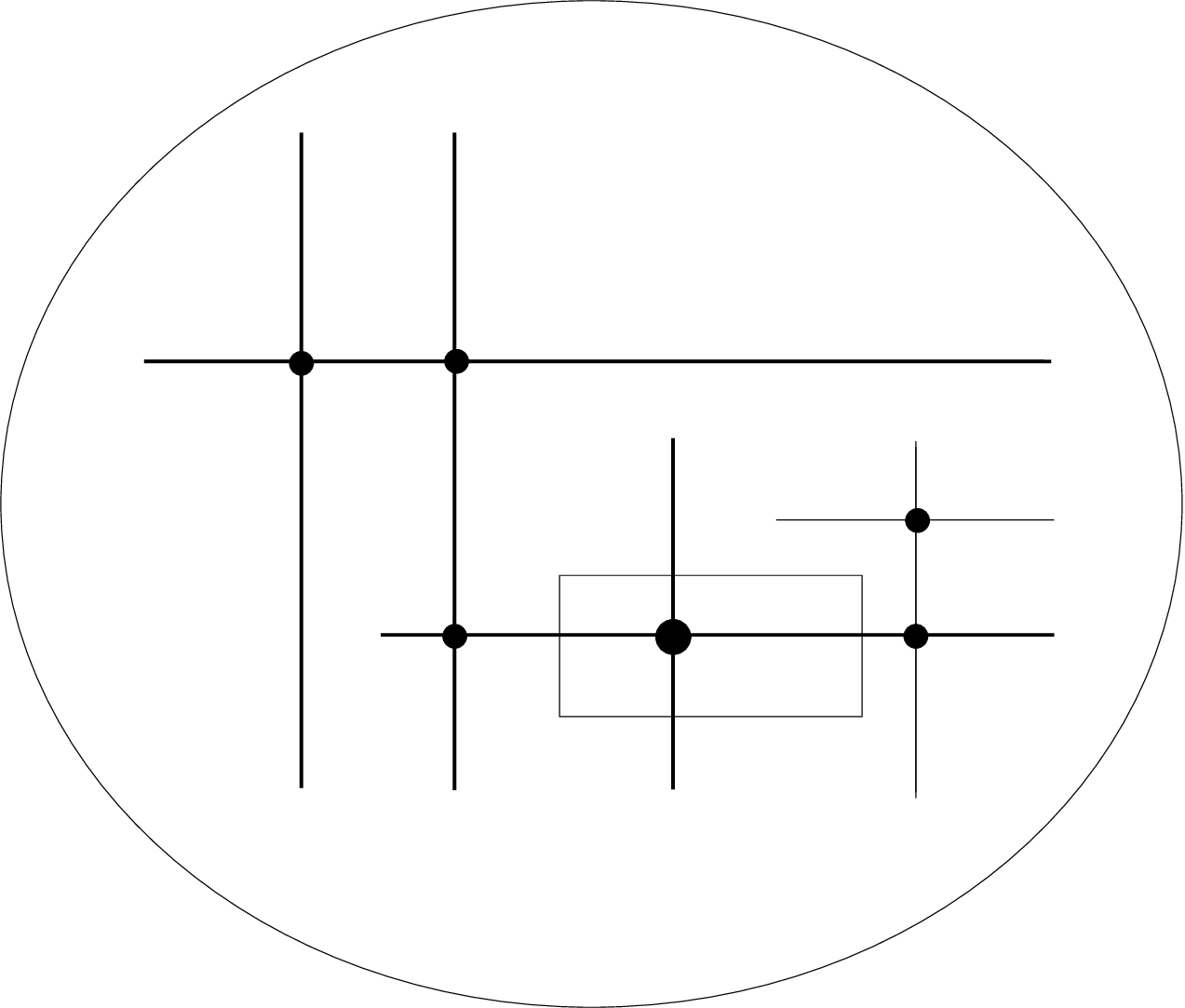}
\end{center}
\caption{The final picture after a resolution of singularities in
the neighborhood of non-normality locus.} \label{figure:sextic solid 7}
\end{figure}

\begin{figure}[htbp]
\begin{center}
\includegraphics[width=.4\textwidth]{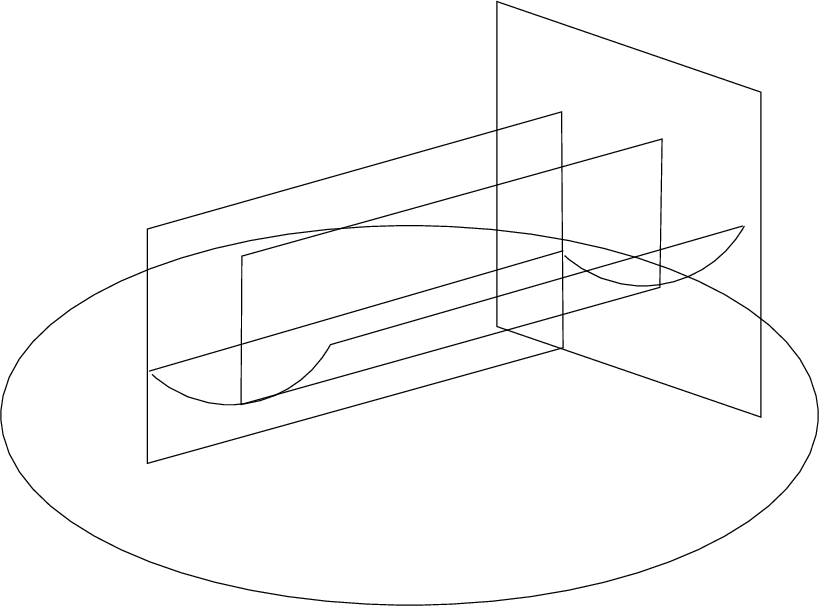}
\end{center}
\caption{The resolution in the neighborhood of the ``first sheet'' of Figure~\ref{figure:sextic solid vertical}.} \label{figure:sextic solid 8}
\end{figure}

\begin{figure}[htbp]
\begin{center}
\includegraphics[width=.5\textwidth]{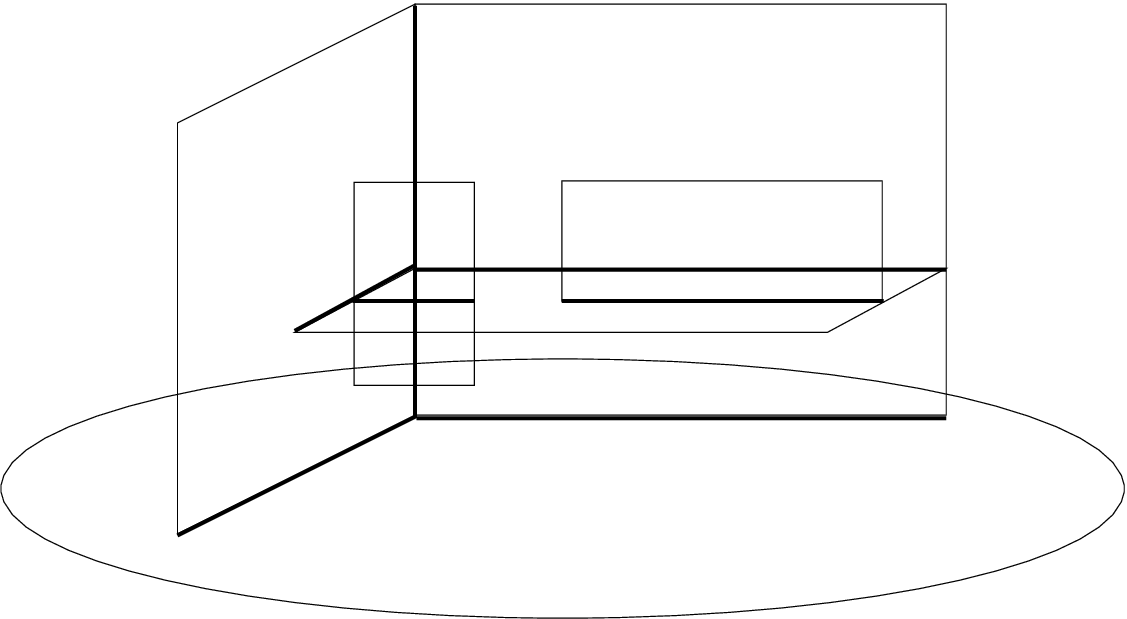}
\end{center}
\caption{The resolution in the neighborhood of the ``deep sheet'' of Figure~\ref{figure:sextic solid vertical}.}
\label{figure:sextic solid 9}
\end{figure}

\begin{figure}[htbp]
\begin{center}
\includegraphics[width=.3\textwidth]{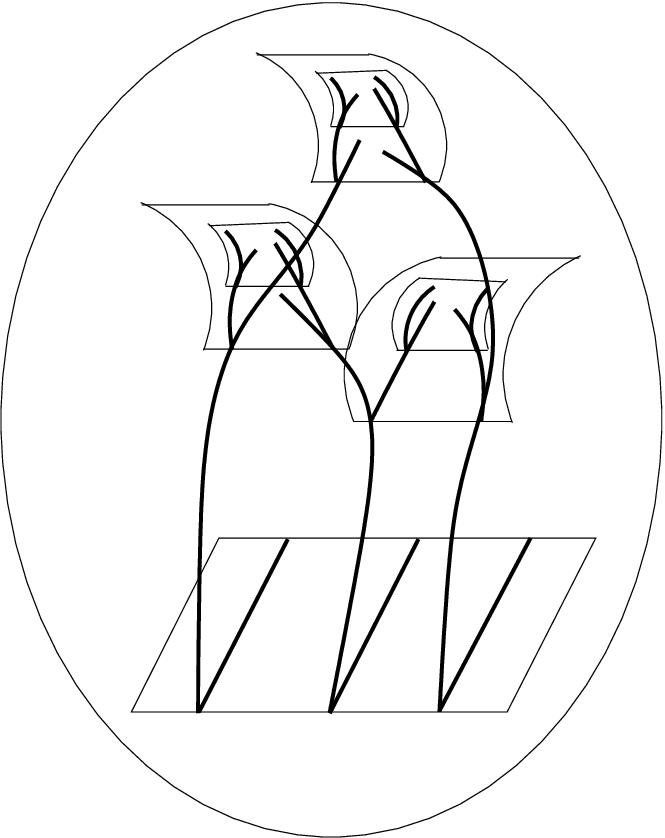}
\end{center}
\caption{After blowing up of horizontal singularities.}
\label{figure:sextic solid vertical}
\end{figure}

We normalize the  plane of singularities blowing it up (twice). Then we  resolve horizontal  vertical singularities.
We record  the  structure of  central fiber and vertical singularities on Figures~\ref{figure:sextic solid 7},~\ref{figure:sextic solid 9},
and~\ref{figure:sextic solid 8} glued in a way given by Figure~\ref{figure:sextic solid vertical}.

The lines on Figure~\ref{figure:sextic solid 7} are surfaces (we look on them ``from above'').
Bold ones intersect the ``base'' surface.
The rectangle is a surface lying ``over'' the ``base''.
It intersects in  two curves (which do not
intersect the base surface) two remaining   surfaces. The point of intersection of these lines and a
``vertical'' line of intersection of two other planes is
denoted by a  fat point. At the end  we get  nine surfaces and twelve lines
recorded on the picture.

We follow procedure of resolving singularities as in previous examples.
The final picture is  obtained  by gluing configurations of surfaces drawn on Figures~\ref{figure:sextic solid 7},~\ref{figure:sextic solid 9},
and~\ref{figure:sextic solid 8} along Figure~\ref{figure:sextic solid vertical}.
More detailed description of Landau--Ginzburg model for sextic double solid see in~\cite{MARCOR}.  Direct calculations   (\cite{MARNON}, \cite{KNS})  produce.

\begin{proposition}  The monodromy of the singular fiber at zero of the  Landau--Ginzburg model for a sextic
double solid with 35 singular points is strictly unipotent.
\end{proposition}

\medskip

The results from \cite{MARNON} suggest that  double covering of quadric ramified in octic with 20 nodal singular points will also have strictly unipotent  monodromy of the singular fiber at zero of its   Landau--Ginzburg model. Indeed this double covering is  nothing else but a three dimensional  quartic deformation and its monodromy was computed in   \cite{MARNON}.
We extract categorical information from this common phenomenon --- strict unipotency of monodromy in  following theorems and  conjectures.

Let us denote by $H(LG(X), \mathcal{F})$ the hypercohomologies of the
perverse sheaf of vanishing cycles on the Landau--Ginzburg model.
 $H(LG(X), \mathcal{F})$ measure cohomologies of $X$ and the monodromy of
 $LG(X)$ --- see \cite{GK} and \cite{KRIS1}.


\begin{theorem} Let $X$ be a smooth Fano variety. Let $LG(X)$ be its Landau--Ginzburg model (in particular, HMS for $X$ and $LG(X)$ holds).
Then the Hochschild homology of Fukaya--Seidel category of $LG(X)$
is $H(LG(X), \mathcal{F}) $.
\end{theorem}

\begin{proof}
It follows from \cite{KKP1}.
\end{proof}

According to Homological Mirror Symmetry the Hochschild homology of Fukaya--Seidel category of Landau--Ginzburg model are isomorphic to   Hochschild homology of  category $D^b(X)$.

Combining  results from Subsections~\ref{LG for quartic2solid} and~\ref{LG for V10} with  conifold transition change  described in \cite{KRIS1} we get

\begin{proposition}  The  Hochschild homology of  category $D^b(X)$ of   Artin--Mumford example and of  resolved $V_{10}$ with  10 singular points are isomorphic.
\end{proposition}

\begin{proof} This follows from a direct calculations of the cohomology
of resolved $V_{10}$ with  10 singular points.
\end{proof}

In fact these   homology  look like cohomologies of a projective space.

Using above  analysis of monodromy  of  Landau--Ginzburg models of
  Artin--Mumford example, $V_{10}$ with  10 singular points, of
 double covering of quadric ramified in octic with 20 nodal singular points,
and of  double solid with ramification in a sextic with 35 singular points
(see \cite{MARNON}) we arrive at

\begin{conjecture}  \label{conjecture:Enriques} The  categories  $D^b(X)$ of the Artin--Mumford example, of   $V_{10}$ with  10 singular points, of
 double covering of quadric ramified in octic with 20 nodal singular points,
and of  double solid with ramification in a sextic with 35 singular points
 contain  category of a nodal
Enriques surface as a semi-orthogonal summand.
\end{conjecture}

\begin{remark} While this paper was being written
Ingalls and Kuznetsov, familiar with our work, stated above conjecture for  Artin--Mumford example
and proved it for the minimal resolution of this example
--- \cite{KC}. The first two authors are collaborating with A.\,Kuznetsov in order to prove this conjecture for  $V_{10}$ with 10 singular points.

\end{remark}

In the next section we look at the above observations from prospective of  theory of spectra   of category.

\section{Spectrum, enhanced spectrum and applications}

\subsection{Classical Spectrum}
In this subsection we review the notions of spectra and gaps following~\cite{BFK}.

Noncommutative Hodge structures were introduced in \cite{KKP1}, as a
means of bringing the techniques and tools of Hodge theory into the
categorical and noncommutative realm.
In the classical setting much of the information about an isolated
singularity is recorded by means of the Hodge spectrum, a set of
rational eigenvalues of the monodromy operator.
A categorical analogue of this Hodge spectrum appears in the
works of Orlov and Rouqier \cite{OR}, \cite{RU}. Let us call this \emph{the Orlov spectrum}.
 Recent work in the manuscript \cite{BFK},
suggests an intimate connection with the classical singularity  theory.

Let us recall the definitions of the Orlov spectrum and discuss some of
the main results in \cite{BFK}.  Let $\mathcal T$ be a triangulated
category.
For any $ G \in \mathcal T $ denote by $ \langle G \rangle_0 $
the smallest full subcategory containing $ G $ which is closed under
isomorphisms, shifting, and taking finite direct sums and summands.
Now inductively define $ \langle G \rangle_n $ as the full subcategory of
objects, $ B $, such that there is a distinguished triangle,
$ X \to B \to Y \to X[1] $, with $ X \in \langle G \rangle_{n-1} $ and
$ Y \in \langle G \rangle_0 $, and direct summands of such objects.

\begin{definition}
Let $G$ be an object of a triangulated category $\mathcal{T}$.  If
there is an $n$ with $\langle G \rangle_{n} = \mathcal T$, we set
\begin{displaymath}
 t(G)=  \text{min } \lbrace  n \geq 0 \  | \ \langle G
 \rangle_{n} = \mathcal T \rbrace.
\end{displaymath}
Otherwise we set $t(G) = \infty$.  We call $t(G)$ the
\emph{generation time} of $G$. If $t(G)$ is finite, we say
that $G$ is a \emph{strong generator}. The \emph{Orlov spectrum}
of $\mathcal T$ is the union of all possible generation times for
strong generators of $\mathcal T$.  The \text{Rouqier dimension} is
the smallest number in the Orlov spectrum.  We say that a triangulated
category $\mathcal T$ has a \emph{gap} of length $s$ if $a$ and
$a+s$ are in the Orlov spectrum but $r$ is not in the Orlov spectrum
for $a < r < a+s$. We denote the maximum (finite) gap of the Orlov spectrum of $\mathcal T$
by $\Gap(\mathcal T)$.
\end{definition}

The following 3 conjectures are from~\cite{BFK}.

\begin{conjecture} \label{gap bound}
If $X$ is a smooth variety then any gap of $D^{b}(X)$ is at most
the Krull dimension of $X$.
\end{conjecture}


\begin{conjecture} \label{conjecture: BFK1} The maximal gap in Orlov's spectrum  is a  birational invariant.
\end{conjecture}


In particular, this conjecture says that if $X$ is smooth projective rational threefold   then  gap of $D^{b}(X)$ is equal to 1.

We now apply theory of gaps to the observations from the previous sections.
First we formulate:

\begin{conjecture} Let $X$ be a smooth algebraic surface. Then $h^{2,0}(X)=0$
is equivalent to $\Gap(D^b(X))=1$.
\end{conjecture}

Combining this  conjecture with Conjecture~\ref{conjecture:Enriques}
we get

\begin{conjecture}   The  gap of the  category $D^b(X)$ for the Artin--Mumford example, of $V_{10}$ with  10 singular points,
 of the double covering of quadric ramified in octic with 20 nodal singular points, and of the double solid with ramification in a sextic with 35 singular points is  equal to 1.
\end{conjecture}

In other words the gap of Orlov spectra is too weak of a categorical  invariant to distinguish the rationality of these  examples. In the next section we introduce  more advanced Noether--Lefschetz spectra.

\subsection{ Enhanced  Noether--Lefschetz Spectra}

Let $\mathcal T$ be an enhanced triangulated category and let $HH^*(\mathcal T)$ be its Hochschild cohomology.

\begin{figure}[h]
  \begin{center}
\includegraphics[width=1.8\nanowidth]{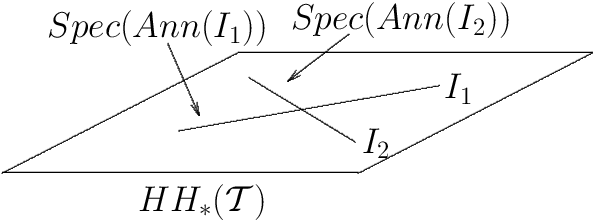}
\caption{Noether--Lefschetz spectra.}\label{fig:2MP}
  \end{center}
\end{figure}

\begin{definition} We denote by \emph{Noether--Lefschetz spectra} $NL(\mathcal T)$ the ordered  collection of sets   over    $HH^*(\mathcal T)$  defined as follows.
For any graded ideal $I$ in   $HH^*(\mathcal T)$ we consider the DG subcategory $Ann(I)$ in $\mathcal T$  --- the annihilator of $I$. The set $\Spec(Ann(I))$ is the set of generators of $\mathcal T$ in  the DG subcategory $Ann(I)$. We
denote the  maximum gap of $\Spec(Ann(I))$ over all subsets  $I$  by $\NLGap(\mathcal T)$ (see Figure~\ref{fig:2MP}).
\end{definition}

Clearly  $\Spec(\mathcal T)$ embeds in the set $(I, \Spec(Ann(I)))$ but the behavior of the gaps in  $NL(\mathcal T)$ is much more complex. For more examples see
\cite{FK}.

We make  the following:

\begin{conjecture} Let $X$ be a 3-dimensional smooth projective variety.
If $X$ is rational then the  gaps in $NL(D^b(X))$ are equal to 1.
\end{conjecture}

The above conjecture suggests  a new invariant of rationality. It is based
on our studies of Landau--Ginzburg models from  previous sections. Theorem 5.1 together with HMS suggests that  $NL(D^b(X))$
are completely determined by monodromy and vanishing cycles of Landau--Ginzburg models, see Table~\ref{tab:Enspec11}.
 Still it is possible that   $NL(D^b(X))$ has all gaps equal to one and $X$ is not rational.

\begin{table}[h]
\begin{center}
\begin{tabular}{|c|c|}
\hline
%
%
%
\begin{minipage}[c]{1in}
\centering
\medskip

A category $\mathcal T$

\medskip

\end{minipage}
&
\begin{minipage}[c]{4in}
\centering
\medskip

$NLSpec (\mathcal T)$

\medskip

\end{minipage}
\\\hline\hline
\begin{minipage}[c]{1in}
\centering
\medskip

$D^b(X)$

\medskip

\end{minipage}
&
\begin{minipage}[c]{4in}
\centering
\medskip

$NLSpec (\mathcal T)\subset Spec_{dg-gr}(HH^*(D^b(X)))\times \Spec (\mathcal T)$

\medskip

\end{minipage}
\\\hline
%
%
%
%
%
%
%
%
%
%
%
%
\begin{minipage}[c]{1in}
\centering
\medskip

$FS(LG(X))$

\medskip

\end{minipage}
&
\begin{minipage}[c]{4in}
\centering
\medskip

$NLSpec (\mathcal T)\subset Spec_{dg-gr}(H^*(LG(X),\mathcal F))\times \Spec (\mathcal T)$

\medskip

\end{minipage}
\\\hline
\end{tabular}
\end{center}
\caption{ HMS and Noether--Lefschetz spectra.}
\label{tab:Enspec11}
\end{table}

In what follows we give conjectural examples of 3-dimensional varieties which have gaps equal to 1
in $\Spec(D^b(X))$ and  have gaps  equal to two or  higher  in $NL(D^b(X))$.
Following  Conjecture~\ref{conjecture: BFK1} 
Homological Mirror Symmetry, and examples in Section 5 we make

\begin{conjecture} In all examples: Artin--Mumford example,   $V_{10}$ with  10 singular points,  double covering of quadric ramified in octic with 20 nodal singular points,
and  double solid with ramification in a sextic with 35 singular points
$\NLGap (D^{b}(X))$ is   equal to two or  higher.
\end{conjecture}

This conjecture is based on the  fact
that Landau--Ginzburg models for
  Artin--Mumford example, for   $V_{10}$ with  10 singular points, for
 double covering of quadric ramified in octic with 20 nodal singular points
and for  double solid with ramification in a sextic with 35 singular points
 have the same monodromies --- see also \cite{MARNON}.

We record all our findings and  conjectures in Table~\ref{tab:ENSPEC22}.

\begin{table}[h]
\begin{center}
\begin{tabular}{|c|c|c|c|}
\hline
\begin{minipage}[c]{2in}
\centering
\medskip

A Fano variety $X$

\medskip

\end{minipage}
&
\begin{minipage}[c]{1.5in}
\centering
\medskip

$D^b(X)$ and $HH_0(X)$

\medskip

\end{minipage}
&
\begin{minipage}[c]{1in}
\centering
\medskip

$\Gap$ $(D^b(X))$

\medskip

\end{minipage}
&
\begin{minipage}[c]{1in}
\centering
\medskip

$\NLGap\,(X)$

\medskip

\end{minipage}
\\\hline\hline
\begin{minipage}[c]{2in}
\centering
\medskip


A double covering of $\PP^3$ ramified in K3 surface with 10 nodal singular points (Artin--Mumford variety).

\medskip

\end{minipage}
&
\begin{minipage}[c]{1.5in}
\centering
\medskip

$D^b(X)=$ $\langle D^b(E), E_1,\ldots,E_{10}\rangle$, where
$E$ is a nodal Enriques surface.


$\dim (HH_0(X))=4$

\medskip

\end{minipage}
&
\begin{minipage}[c]{1in}
\centering
\medskip

$1$

\medskip

\end{minipage}
&
\begin{minipage}[c]{1in}
\centering
\medskip

$\geq 2$

\medskip

\end{minipage}
\\\hline
\begin{minipage}[c]{2in}
\centering
\medskip

Double covering

$\xymatrix{V_{10}\ar[d]^-{\mbox{2:1}}\\V_5}$

\medskip

\end{minipage}
&
\begin{minipage}[c]{1.5in}
\centering
\medskip

$D^b(X)=\langle D^b(E),\ldots\rangle$, where
$E$ is a nodal Enriques surface.


$\dim (HH_0(X))=4$

\medskip

\end{minipage}
&
\begin{minipage}[c]{1in}
\centering
\medskip

$1$

\medskip

\end{minipage}
&
\begin{minipage}[c]{1in}
\centering
\medskip

$\geq 2$

\medskip

\end{minipage}
\\\hline
\begin{minipage}[c]{2in}
\centering
\medskip

$\PP^3$

\medskip

\end{minipage}
&
\begin{minipage}[c]{1.5in}
\centering
\medskip


$\dim (HH_0(X))=4$

\medskip

\end{minipage}
&
\begin{minipage}[c]{1in}
\centering
\medskip

$1$

\medskip

\end{minipage}
&
\begin{minipage}[c]{1in}
\centering
\medskip

$1$

\medskip

\end{minipage}
\\\hline
\begin{minipage}[c]{2in}
\centering
\medskip

Sextic double solid with 35 nodal singular points.

\medskip

\end{minipage}
&
\begin{minipage}[c]{1.5in}
\centering
\medskip


$D^b(X)=\langle D^b(E),\ldots\rangle$, where
$E$ is a nodal Enriques surface.

$\dim (HH_0(X))=4$

\medskip

\end{minipage}
&
\begin{minipage}[c]{1in}
\centering
\medskip

$1$

\medskip

\end{minipage}
&
\begin{minipage}[c]{1in}
\centering
\medskip

$\geq 2$

\medskip

\end{minipage}
\\\hline
\begin{minipage}[c]{2in}
\centering
\medskip

Double covering of quadric ramified in octic with 20 nodal singular points.

\medskip

\end{minipage}
&
\begin{minipage}[c]{1.5in}
\centering
\medskip


$D^b(X)=\langle D^b(E),\ldots\rangle$, where
$E$ is a smooth Enriques surface.

$\dim (HH_0(X))=4$

\medskip

\end{minipage}
&
\begin{minipage}[c]{1in}
\centering
\medskip

$1$

\medskip

\end{minipage}
&
\begin{minipage}[c]{1in}
\centering
\medskip

$\geq 2$

\medskip

\end{minipage}
\\\hline
\end{tabular}
\end{center}
\caption{Summarizing conjectures.}
\label{tab:ENSPEC22}
\end{table}

\begin{remark}
\label{remark:AMV10}
It is quite possible that derived categories of the Artin--Mumford example and of $V_{10}$ are related via deformation in which case  equalities  of spectra is not surprising.
\end{remark}

\begin{remark} The considerations in the last two sections suggest a strong correlation between spectra, monodromy and walls in  moduli spaces of stability conditions. We pose the following two questions:

\medskip

{\bf Question 1.}  Does Noether--Lefschetz spectra define a stratification on
the moduli space of stability conditions?

\medskip

{\bf Question 2.}  Are classical  Noether--Lefschetz loci  connected to this  stratification?

\end{remark}

\begin{remark} Artin--Mumford example is an example of a conic bundle.
We expect that technique discussed here will lead to many examples of conic bundles for which the gap of Orlov's spectrum is equal to  one and their nonrationality can be established using  gaps in Noether--Lefschetz spectra.

\end{remark}


\bigskip

\address{
Atanas Iliev, SNU,
}{email: ailiev2001@yahoo.com}

\medskip

\address{
Ludmil Katzarkov,
UW,
}{email: lkatzark@math.uci.edu}

\medskip

\address{
Victor Przyjalkowski, MI RAS,
}{email: victorprz@mi.ras.ru, victorprz@gmail.com}

\end{document}